\documentclass[10pt]{agn_article}
\usepackage[bibtex]{agn_all}
\usepackage[multiple]{footmisc}
\usepackage{longtable}
\usepackage{makecell}
\usepackage{marvosym}
\usepackage{multicol}
\usepackage{tikz}
\usepackage{amsmath}
\usepackage{agn_tikz}
\usepackage{caption}
\def\sample{500}
\newcommand{\blue}[1]{{\color{blue}#1}}
\newcommand{\red}[1]{{\color{red}#1}}
\renewcommand{\rho}{\varrho}

\renewcommand{\div}{\mathrm{Div}}
\newcommand{\skw}{\mathrm{skew}}

\newcommand*\dif{\mathop{}\!\mathrm{d}}

\def\C{\mathbb{C}}
\def\sk{\textnormal{skew}}
\DeclareMathOperator{\Biot}{Biot}

\DeclareMathOperator{\ZJ}{ZJ}

\DeclareMathOperator{\bfsym}{\textbf{sym}}
\DeclareMathOperator{\iso}{iso}
\DeclareMathOperator{\NH}{NH}
\DeclareMathOperator{\lin}{lin}
\newcommand{\DD}{\mathrm{D}}
\newcommand{\WW}{\mathrm{W}}

\newcommand{\notiff}{%
	\mathrel{{\ooalign{\hidewidth$\not\phantom{"}$\hidewidth\cr$\iff$}}}}
\def\dd{\displaystyle}

\numberwithin{equation}{section}

\makeatletter
\DeclareFontFamily{OMX}{MnSymbolE}{}
\DeclareSymbolFont{MnLargeSymbols}{OMX}{MnSymbolE}{m}{n}
\SetSymbolFont{MnLargeSymbols}{bold}{OMX}{MnSymbolE}{b}{n}
\DeclareFontShape{OMX}{MnSymbolE}{m}{n}{
	<-6>  MnSymbolE5
	<6-7>  MnSymbolE6
	<7-8>  MnSymbolE7
	<8-9>  MnSymbolE8
	<9-10> MnSymbolE9
	<10-12> MnSymbolE10
	<12->   MnSymbolE12
}{}
\DeclareFontShape{OMX}{MnSymbolE}{b}{n}{
	<-6>  MnSymbolE-Bold5
	<6-7>  MnSymbolE-Bold6
	<7-8>  MnSymbolE-Bold7
	<8-9>  MnSymbolE-Bold8
	<9-10> MnSymbolE-Bold9
	<10-12> MnSymbolE-Bold10
	<12->   MnSymbolE-Bold12
}{}

\let\llangle\@undefined
\let\rrangle\@undefined
\DeclareMathDelimiter{\llangle}{\mathopen}%
{MnLargeSymbols}{'164}{MnLargeSymbols}{'164}
\DeclareMathDelimiter{\rrangle}{\mathclose}%
{MnLargeSymbols}{'171}{MnLargeSymbols}{'171}
\makeatother
\title{The corotational stability postulate: positive incremental Cauchy stress moduli for diagonal, homogeneous deformations in isotropic nonlinear elasticity}
\defineaffiliation{tud}{
	Faculty of Mathematics,
	Technical University of Dresden,
	Zellescher Weg 12--14,
	01069 Dresden, Germany
}
\defineaffiliation{Lavrentyev}{
Sergey N. Korobeynikov, Principal Investigator of Composite Mechanics Laboratory of Lavrentyev Institute of Hydrodynamics, Lavrentyev Prospekt 15, Novosibirsk, 630090, Russia, email: S.N.Korobeynikov@mail.ru
}
\defineaffiliation{Maroua}{
Aurélien S. Nguetcho Tchakoutio, Départment of Physics, Faculty of Science, University de Maroua, Po. Box: 814 Maroua - Cameroon, email: nguetchoserge@yahoo.fr
}
\author{%
	Patrizio Neff\thanks{
		Patrizio Neff, University of Duisburg-Essen, Head of Chair for Nonlinear Analysis and Modelling, Faculty of Mathematics, Thea-Leymann-Straße 9, D-45127 Essen, Germany, email: patrizio.neff@uni-due.de
	}
	,\quad
	Nina J. Husemann\thanks{
		Nina J. Husemann, University of Duisburg-Essen, Chair for Nonlinear Analysis and Modelling, Faculty of Mathematics, Thea-Leymann-Straße 9, D-45127 Essen, Germany, email: nina.husemann@stud.uni-due.de
	}
	,\\
	Aurélien S.  Nguetcho Tchakoutio\affil{Maroua}
	, \quad
	Sergey N. Korobeynikov\affil{Lavrentyev}
	\quad and\quad
	Robert J. Martin\thanks{
		Robert J. Martin, University of Duisburg-Essen, Chair for Nonlinear Analysis and Modelling, Faculty of Mathematics, Thea-Leymann-Straße 9, D-45127 Essen, Germany, email: robert.martin@uni-due.de
	}
}
\begin{document}
\maketitle
 \begin{abstract}
\noindent In isotropic nonlinear elasticity the corotational stability postulate (CSP) is the requirement that
\begin{equation*}
	\langle\frac{\DD^{\circ}}{\DD t}[\sigma] , D \rangle > 0 \quad \forall \ D \in \Sym(3)\setminus \{0\} \, ,
\end{equation*}	
where $\frac{\DD^{\circ}}{\DD t}$ is \textbf{any} corotational stress rate, $\sigma$ is the Cauchy stress and $D = \sym L$, condition $L= \dot{F} \, F^{-1}$ is the deformation rate tensor. For $\widehat{\sigma} (\log V) \colonequals \sigma (V)$ it is equivalent to the monotonicity (TSTS-M$^+$)
\begin{equation*}
	\langle \widehat{\sigma} (\log V_1) - \widehat{\sigma} (\log V_2) , \log V_1 - \log V_2 \rangle > 0 \quad \forall \ V_1, V_2 \in \Sym^{++}(3), \ V_1 \neq V_2 \, .
\end{equation*}
For hyperelasticity, (CSP) is in general independent of convexity of the mapping $F \mapsto \WW(F)$ or $U \mapsto \widehat{\WW}(U)$.
Considering a family of diagonal, homogeneous deformations $t \mapsto F(t)$
one can, nevertheless, show that (CSP) implies positive incremental Cauchy stress moduli for this deformation family, including the incremental Young's modulus, the incremental equibiaxial modulus, the incremental planar tension modulus and the incremental bulk modulus. Aside, (CSP) is sufficient for the Baker-Ericksen and tension-extension inequality. Moreover, it implies local invertibility of the Cauchy stress-stretch relation. Together, this shows that (CSP) is a reasonable constitutive stability postulate in nonlinear elasticity, complementing local material stability viz. LH-ellipticity.
\end{abstract}
\vspace*{.49em}%
\vspace*{.49em}%
\clearpage
\tableofcontents

\section{Introduction}
Recently, the corotational stability postulate (CSP) has been introduced into nonlinear elasticity \cite{CSP2024}.
It amounts to a constitutive requirement expressed in rate-type format, namely
\begin{equation}
	\label{eq:CSP}
	\langle \frac{\DD^{\circ}}{\DD t} [\sigma] , D \rangle > 0 \quad \forall \ D \in \Sym(3)\setminus \{0\}
\end{equation}
where $\frac{\DD^{\circ}}{\DD t}$ is \textbf{any} reasonable corotational rate (cf.~\cite{poscor2024})
\begin{equation}
	\frac{\DD^{\circ}}{\DD t} [\sigma] = \frac{\DD} {\DD t}[\sigma] + \sigma \, \Omega^{\circ} - \Omega^{\circ} \, \sigma  \, , \quad \Omega^{\circ} \in \mathfrak{so}(3) \, ,
\end{equation}
$\sigma$ is the Cauchy stress tensor and $D=\sym L$ is the deformation rate tensor, where $L=\dot{F} \, F^{-1}= \DD_{\xi}v$ is the spatial velocity gradient and $\Omega^{\circ}$ is a spin-tensor. For more information on the notation please consult the Appendix \ref{sec:app_notation}.\\

In \cite{CSP2024} it has been shown that \eqref{eq:CSP} is equivalent\footnote{
for $\frac{\DD^{\circ}}{\DD t}= \{ \frac{\DD^{\ZJ}}{\DD t}, \frac{\DD^{\log}}{\DD t}\}$. Here, $\frac{\DD^{\ZJ}}{\DD t} [\sigma] = \frac{\DD}{\DD t} [\sigma] + \sigma \, W - W \sigma, \, W = \skw \, \text{L}$ is the Zaremba-Jaumann rate. The result for all reasonable corotational rates is about to be submitted (cf.~\cite{Martin2024}).
}
to the already introduced (TSTS-M$^{++}$) condition (cf. \cite{agn_neff2015exponentiatedI,agn_neff2015exponentiatedII,jog2013}) which reads for $\widehat{\sigma}(\log V) \colonequals \sigma (V)$
\begin{equation}
\label{eq:TSTS-M++}
	\text{TSTS-M$^{++}$:} \quad \sym \, \DD_{\log V} \widehat{\sigma} (\log V) \in \Sym_4^{++}(6) \, .
\end{equation}
Note that $\DD_{\log V} \widehat{\sigma} (\log V)$ is not necessarily major-symmetric (cf.~\cite{majorsym2024}).
The latter implies further the Hilbert-monotonicity
\begin{equation}
	\label{eq:TSTS-M+}
	\text{TSTS-M$^{+}$:} \quad \langle \widehat{\sigma} (\log V_1) - \widehat{\sigma} (\log V_2) , \log V_1 - \log V_2 \rangle > 0 \quad \forall \ V_1, V_2 \in \Sym^{++}(3), \ V_1 \neq V_2
\end{equation}
so that \eqref{eq:TSTS-M+} is one possibility to express that ''stress increases with strain'' in fact: ``Cauchy stress $\sigma$ increases with logarithmic strain``. It is already known that (CSP) implies the BE-inequalities (cf. \cite{bakerEri54}) and the tension-extension TE-inequalities (cf. \cite{leblond1992constitutive,Leblond2024}). Moreover, TSTS-M$^{++}$ is easily checked in the isotropic case by switching to the representation in principal stresses versus principal logarithmic stretches, see \cite{Ghiba2024,Martin2024Hill}.\\

For incompressibility, condition \eqref{eq:TSTS-M+} turns into the well-known Hill's inequality (cf. \cite{CSP2024,sidoroff1974restrictions,hill1957,hill1968constitutive,hill1970constitutive}) for the Kirchhoff-stress $\tau= \det F \cdot \sigma$, and $\widehat{\tau} (\log V) \colonequals \tau (V)$
\begin{equation}
	\label{eq:taumonoton}
	\langle \widehat{\tau}(\log V_1) - \widehat{\tau} (\log V_2) , \log V_1 - \log V_2 \rangle > 0 \quad \forall \ V_1, V_2 \in \Sym^{++}(3), \ V_1 \neq V_2 \cdot
\end{equation}
Note that \eqref{eq:taumonoton} is, similarly  to \eqref{eq:CSP}, equivalent to the rate-condition
\begin{equation}
	\label{eq:Hillsinequality}
	\langle \frac{\DD^{\ZJ}}{\DD t} [\tau] , D \rangle > 0 \quad \forall \ D \in \Sym(3) \setminus \{0\} \, ,
\end{equation}
where $\frac{\DD^{\ZJ}}{\DD t}$ denotes the Zaremba-Jaumann corotational rate (cf.~\cite{Leblond2024}).\\

A linearized version of \eqref{eq:CSP} for small stress and small rotations is the requirement of ''positive second order internal work'':
\begin{equation}
	\label{eq:positivesecondorderwork}
	\langle \dot{\sigma}^{\lin} , \dot{\varepsilon} \rangle > 0 \quad \forall \ \dot{\varepsilon} \in \Sym(3)\setminus \{0\}
\end{equation}
which is equivalent to
\begin{equation}
	\label{eq:CisoSym}
	\C^{\iso} \in \Sym_4^{++}(6)
\end{equation}
for the linear elastic constitutive law
\begin{equation}
	\sigma^{\lin}= \C^{\iso}. \varepsilon = 2 \, \mu \, \varepsilon + \lambda \, \tr (\varepsilon) \, \id \, , \quad \quad \mu > 0 , \ 2 \, \mu + 3 \, \lambda > 0 \, .
\end{equation}
While the meaning of \eqref{eq:CisoSym} is established as the stability requirement for linear elasticity and \eqref{eq:positivesecondorderwork} implies that stress increases with strain in physically nonlinear, but geometrically linear problems (as e.g. small strain plasticity, cf. Figure \ref{fig:infDruckerstability}), the immediate interpretation of conditions \eqref{eq:CSP}, \eqref{eq:TSTS-M++}, \eqref{eq:TSTS-M+} is less clear, owing to the appearance of the (arbitrary) corotational rate $\frac{\DD^{\circ}}{\DD t}$ and being defined in the spatial configuration.\\

Here, we will therefore simplify the setting by considering a family of deformations $t \mapsto \varphi (t)$ of the domain $V_0 \subset \R^3$ such that $t \mapsto F(t) =\DD \varphi (t) = \diag (\lambda_1 (t) , \lambda_2 (t), \lambda_3 (t))$ is diagonal, and therefore has \textit{constant principal axes}.
Furthermore, we will assume that $F(t)$ is homogeneous in the space variable (thus a subset of so called universal deformations, cf. \cite{yavari2024,Ericksen1955}) and we assume hyperelastic response.
Let us remark that the assumed loading deformation history corresponds to the standard homogeneous test protocolls, like uniaxial tension, planar tension, equibiaxial tension etc. (but not simple shear as for simple shear, $F$ is not diagonal and the principal axes rotate).\\
Our main finding is that for the admitted tests, (CSP) implies a priori positive incremental Cauchy-stress moduli. 
In \cite{Young2024} it is shown that (CSP) also implies a positive incremental pure shear modulus. Whether this remains true for (not so simple) simple shear \cite{thiel2018} will be investigated in another contribution.\\
Note carefully that the above conclusion is only true in the assumed circumstances while in general (CSP) does neither imply convexity of $F \mapsto \WW(F)$ nor convexity of $U \mapsto \widehat{\WW}(U)$ as the example of the exponentiated Hencky energy (cf. \cite{agn_neff2015exponentiatedI})
\begin{equation}
	\WW_{\text{exp-Hencky}} (F)= \frac{\mu}{k} \exp\left(k \, \norm{\log V}^2\right) + \frac{\lambda}{2 \, \widehat{k}} \exp\left(\widehat{k} \, (\log(\det V))^2\right)
\end{equation}
shows: $\WW_{\text{exp-Hencky}}$ satisfies (CSP) throughout but is, of course, not convex in $F$ and indeed not convex in $U= \sqrt{F^T \, F}$.
It is also known that (CSP) is otherwise independent of LH-ellipticity and polyconvexity (cf. \cite{CSP2024,leblond1992constitutive,Leblond2024}) since $\WW_{\text{exp-Hencky}}$ is not LH-elliptic everywhere. It should also be observed that a standard compressible Neo-Hooke model does not satisfy (CSP) while the incompressible Neo-Hooke model complies with (CSP) and is polyconvex.\\
\section{Second order internal work condition versus corotational stability in the spatial and referential picture} \label{appsecondorderwork}
In this section we will have a closer look at the concept of ``second order internal work'' since it bears some superficial  resemblance to the (CSP) condition. It proves useful, however, to start with linear elasticity.
\subsection{Linear elasticity}
In linear elasticity, the stored elastic energy can be expressed as
\begin{align}
	\mathcal{E}^{\lin}(t) = \int_{V_0} \WW^{\lin}(\varepsilon(t)) \dif V_0 = \int_{V_0} \frac12 \, \langle \C . \varepsilon(t), \varepsilon(t) \rangle \dif V_0 = \int_{V_0} \frac{1}{2} \langle \sigma^{\lin} (t) , \varepsilon (t) \rangle \dif V_0 \qquad \textnormal{``work''}.
\end{align}
Let us consider the expansion
\begin{equation}
	\begin{alignedat}{3}
		\sigma^{\lin} (t + \delta t) &= \sigma^{\lin} (t) +  \dot{\sigma}^{\lin}(t) \, \delta t + \text{h.o.t.} &&= \sigma^{\lin} (t) + \dif \sigma^{\lin}(t) + \text{h.o.t.} \\
		\varepsilon (t + \delta t) &= \varepsilon (t) + \dot{\varepsilon} (t) \, \delta t + \text{h.o.t.} &&=\varepsilon (t) + \dif \varepsilon(t) + \text{h.o.t.} \, .
	\end{alignedat}
\end{equation}
According to Petryk\footnote{
	In \cite[p.377f]{petryk1985}:
		``The second-order work of deformation is a classical concept in the theory of plasticity. In the so-called small strain theory, or more precisely, when geometry changes are disregarded, the second-order work per unit volume during {\em proportional} application of a small increment $\delta \varepsilon_{ij}$ in strain components is, by definition, equal to $\frac{1}{2} \delta \sigma_{ij} \, \delta \varepsilon_{ij}$ (or $\frac{1}{2} \delta \sigma \cdot \delta \varepsilon$ in the symbolic notation), where $\delta \sigma_{ij}$ are the respective small increments of the stress components; the summation convention is used for repeated subscripts. That expression plays a fundamental role in Drucker's \cite{Drucker1950,Drucker1951} definition of work-hardening, interpreted as a postulate of stability of the material in a restricted sense. A similar expression, integrated over the body volume, appears in Hill's \cite{Hill1958,Hill1959} condition for stability of equilibrium of an inelastic continuous body under dead loading, with geometry changes taken into account.``
} \cite{petryk1985}, the ``second order internal work`` expression in linear elasticity is given by $\frac{1}{2} \langle \dif \sigma^{\lin}, \dif \varepsilon \rangle$ since by expansion for a small strain increment of the elastic energy we have
\begin{equation}
	\begin{alignedat}{2}
		\label{eq:scalarproduct}
		\frac{1}{2} \langle \sigma^{\lin} + \dif \sigma^{\lin} , \varepsilon + \dif \varepsilon \rangle &= \frac{1}{2} \big( \langle \sigma^{\lin} , \varepsilon \rangle + \langle \sigma^{\lin} , \dif \varepsilon \rangle + \langle \dif \sigma^{\lin} , \varepsilon \rangle + \langle \dif \sigma^{\lin} , \dif \varepsilon \rangle\big)\\
		&= \frac{1}{2} \langle \sigma^{\lin} , \varepsilon \rangle + \frac{1}{2} \langle \sigma^{\lin} , \dif \varepsilon \rangle + \frac{1}{2} \langle \dif \sigma^{\lin} , \varepsilon \rangle + \frac{1}{2} \langle \dif \sigma^{\lin} , \dif \varepsilon \rangle \, .
	\end{alignedat}
\end{equation}
Consider similarly the direct expansion of the elastic energy
\begin{equation}
	\label{eq:Elin}
	\mathcal{E}^{\lin} (t + \delta t) = \mathcal{E}^{\lin}(t) + \frac{\dif}{\dif t} \mathcal{E}^{\lin}(t) \, \delta t + \underbrace{\frac{1}{2} \frac{\dif^2}{\dif t^2} \mathcal{E}^{\lin}(t) \, \delta t^2}_{\text{``second order internal work``}} + \, \text{h.o.t.} \, .
\end{equation}
Equating like powers in \eqref{eq:scalarproduct} and \eqref{eq:Elin} suggests already
\begin{equation}
	\frac{1}{2} \frac{\dif^2}{\dif t^2} \mathcal{E}^{\lin}(t) \delta t^2 = \frac{1}{2} \langle \dif \sigma^{\lin} , \dif \varepsilon \rangle \, .
\end{equation}
We will make this idea precise now.
\begin{figure}
	\begin{center}
		\begin{minipage}[h!]{0.45\linewidth}
			\begin{tikzpicture}
				\begin{axis}[scale=1.2,
					axis x line=middle,axis y line=middle,
					width=1\linewidth,
					height=0.6\linewidth,
					xlabel={$\varepsilon$},        
					ylabel={$\sigma$},        
					xtick={-3},
					ytick={-0.2},
					xmin=-2,
					xmax=6,
					ymin=-0.05
					]
					\addplot[blue, thick, domain=0:5.5, samples=\sample]{0.1*ln(x+1)};
					\addplot[red, thick, domain=2:3.5, samples=\sample]{0.0333*x + 0.0431};
					\draw[red,dashed] (axis cs:2,0.1098) -- (axis cs:3.5,0.1098);
					\draw[blue,dashed] (axis cs: 2,0.1098) -- (axis cs:2,0);
					\node at (axis cs: 2,-0.02){\blue{$\varepsilon (t)$}};
					\draw[blue,dashed] (axis cs: 3.5,0.1098) -- (axis cs: 3.5,0);
					\path[draw=black] (axis cs: 2,0) -- (axis cs: 3.5,0) node[pos=0.5,above=-0.2em] {\blue{$\overbrace{\hspace{3.4em}}^{\dif \varepsilon}$}};
					\node at (axis cs: 3.5,-0.02){\blue{$\varepsilon (t + \delta t)$}};
					\draw[blue,dashed] (axis cs: 2,0.1098) -- (axis cs: 0,0.1098);
					\node at (axis cs: -0.4,0.1098){\blue{$\sigma (t)$}};
					\draw[blue,dashed] (axis cs: 3.5,0.15965) -- (axis cs:0,0.15965);
					\node at (axis cs: -0.8,0.15965){\blue{$\sigma (t + \delta t)$}};
					\path (axis cs:-0.5,0.15965) -- (axis cs:-0.5,0.1098) node[pos=0.51, right=-0.3em]{$\blue{\left.\begin{array}{c}\\[0.5em]\end{array}\right\} \dif\sigma}$};
					\draw[red,dashed] (axis cs:3.5,0.15965) -- (axis cs:3.5,0.1098);
					\node at (axis cs: 2.7,0.09){\red{$\dot{\varepsilon} \cdot \delta t$}};
					\node at (axis cs: 4.2,0.13){\red{$\dot{\sigma} \cdot \delta t$}};
				\end{axis}
			\end{tikzpicture}
		\caption{Infinitesimal Drucker stability: Cauchy stress $\sigma$ increases with infinitesimal strain $\varepsilon$ for geometrically linear but physically nonlinear response, as e.g. in work-hardening small strain plasticity.}
		\label{fig:infDruckerstability}
		\end{minipage}
		\qquad
		\begin{minipage}[h!]{0.45\linewidth}
			\begin{tikzpicture}
				\node[draw] at (13,4){\parbox{5cm}{\begin{equation*}
							\begin{alignedat}{2}
								\langle \dot{\sigma} , \dot{\varepsilon} \rangle > 0 \quad &\iff \\ 
								\langle \dif \sigma , \dif \varepsilon \rangle > 0 \quad &\iff \\
								\langle \sigma (\varepsilon_1) - \sigma (\varepsilon_2) , \varepsilon_1 &- \varepsilon_2 \rangle > 0
							\end{alignedat}
						\end{equation*}}};
				\draw[white] (12,2) -- (9,2);
			\end{tikzpicture}
		\caption{Different equivalent expressions for the infinitesimal Drucker stability.}
		\label{fig:infDruckerstabilitybox}
		\end{minipage}
	\end{center}
\end{figure}
Taking time derivatives along a motion yields
\begin{align}
	\frac{\dif}{\dif t} \mathcal{E}^{\lin}(t) &= \int_{V_0} \langle \DD_{\varepsilon} \WW^{\lin}(\varepsilon(t)), \dot{\varepsilon}(t) \rangle \dif V_0 = \int_{V_0} \langle \sigma^{\lin}(t) , \dot{\varepsilon}(t) \rangle \dif V_0 = \mathcal{P}_{\textnormal{int}}^{\lin} \qquad \textnormal{``internal power''} \notag \, , \\
	\frac{\dif^2}{\dif t^2} \mathcal{E}^{\lin}(t) &= \int_{V_0} \langle \dot{\sigma}^{\lin}(t), \dot{\varepsilon}(t) \rangle + \langle \sigma^{\lin}(t), \ddot{\varepsilon}(t) \rangle \dif V_0 \overset{\sigma^{\lin} \in \Sym(3)}{=} \int_{V_0} \langle \dot{\sigma}^{\lin}(t) , \dot{\varepsilon}(t) \rangle + \langle \sigma^{\lin}(t) , \DD u_{,tt} (t) \rangle \dif V_0 \nonumber \\
	&= \int_{V_0} \langle \dot{\sigma}^{\lin}(t), \dot{\varepsilon}(t) \rangle \dif V_0 - \int_{V_0} \langle \underbrace{\textnormal{Div} \, \sigma^{\lin}(t)}_{\substack{= \; 0 \; \textnormal{in} \\ \textnormal{equilibrium}}} , u_{,tt} \rangle \dif V_0  + \int_{V_0} \div (\sigma^{\lin , T} . u_{,tt}) \dif V_0 \nonumber\\
	&= \int_{V_0} \langle \dot{\sigma}^{\lin}(t) , \dot{\varepsilon}(t) \rangle \dif V_0 + \int_{\partial V_0} \langle \sigma^{\lin} . \overset{\rightarrow}{n} , u_{,tt} \rangle \dif V_0 \\
	&= \int_{V_0} \langle \C . \dot{\varepsilon}(t), \dot{\varepsilon}(t) \rangle \dif V_0 + \int_{\partial V_0} \langle \sigma^{\lin} . \overset{\rightarrow}{n} , u_{,tt} \rangle \dif V_0 \ge c^{+} \, \norm{\dot{\varepsilon}(t)}^2_{V_0} + \int_{\partial V_0} \langle \sigma^{\lin} . \overset{\rightarrow}{n} , u_{,tt} \rangle \dif V_0  \, . \notag
\end{align}
For quasistatic loading (or if $u \vert_{\partial V_0}=0$ in order to only consider internal variations of the body), terms with $u_{,tt}$ at the boundary will be dropped (e.g. if $u_{,t}\big\vert_{\partial V} = $ const.), so that we are left with
\begin{equation}
	\frac{\dif^2}{\dif t^2} \mathcal{E}^{\lin}_{\text{int}} (t) = \int_{V_0} \langle \dot{\sigma}^{\lin} (t) , \dot{\varepsilon} (t) \rangle \dif V_0 = \int_{V_0} \langle \C. \dot{\varepsilon} (t) , \dot{\varepsilon}(t) \rangle \dif V_0 \geq c^+ \norm{\dot{\varepsilon} (t)}^2_{L^2(V_0)} \, , \quad \text{if} \ \C \in \Sym_4^{++}(6) \, .
\end{equation}
Thus integrated positive ``second order internal work'' in linear elasticity (or the infinitesimal ``Drucker stability postulate (DSP)``\cite[eq.~(1)]{Drucker1950})\footnote
{Mandel \cite[p.59]{mandel1966}: ''En principe on doit introduire dans cette formule $[\dot{\sigma}_{ij} \, \dot{\varepsilon}_{ij} \geq 0]$ non pas la vitesse de contraintes $\dot{\sigma}_{hk}$ par rapport à des axes fixes, mais la vitesse de contraintes $\frac{\DD \sigma_{hk}}{\DD t}$ par rapport à des axes animés de la vitesse de rotation $\omega_{ij}$ de l'élément matériel. On a:
\begin{equation}
	\frac{\DD^{\ZJ} \sigma_{ij}}{\DD t} = \dot{\sigma}_{ij} - \omega_{ik} \, \sigma_{kj} - \omega_{jk} \, \sigma_{ki} \, .
\end{equation}
Mais nous supposons la vitesse de rotation $\omega$ et les contraintes $\sigma$ suffisamment faibles pour que $\frac{\DD \sigma_{ij}}{\DD t}$ puisse être remplacé par $\dot{\sigma}_{ij}$.''
Translation and update of notation: ''In principle, we must introduce in this formula $\langle \dot{\sigma} , \dot{\varepsilon} \rangle \geq 0$
not the [material] stress rate $\frac{\DD}{\DD t}[\sigma]$ with respect to fixed axes, but rather the stress rate $\frac{\DD^{\ZJ}}{\DD t}[\sigma]$ with respect to axes moving with the rotation rate $W$ of the material element. We have [the Zaremba-Jaumann rate]:
\begin{equation}
	\frac{\DD^{\ZJ}}{\DD t} [\sigma] \colonequals \frac{\DD}{\DD t}[\sigma] + \sigma W - W \sigma \, .
\end{equation}
However, we assume that the rotation rate $W$ and stresses $\sigma$ are sufficiently small so that $\frac{\DD^{\ZJ}}{\DD t}[\sigma]$ can be replaced by $\frac{\DD}{\DD t}[\sigma]$.'' [our comment: and in addition one needs to replace $\dot{\varepsilon}$ by $D= \sym  \DD_{\xi}v$ the deformation rate tensor et voila: the (CSP) emerges (cf. Mandel \cite[p.59]{mandel1966})]. In Mandels text, however, it remains vague, whether $\sigma$ is the Cauchy stress or the Kirchhoff stress.	
}\footnote
{The ``Drucker stability postulate`` \cite[eq.~(1)]{Drucker1950} (DSP) states that the incremental internal energy can only increase upon additional loading. Effectively, it eliminates the possibility of (absolute) strain softening and implies that small strain Cauchy stress increases for increasing stretch in uniaxial tension, cf. Figure \ref{fig:infDruckerstability}. It must be appreciated that Drucker himself did not introduce any geometrically nonlinear generalization of his stability postulate. In the current literature, Druckers stability postulate for the geometrically nonlinear setting seems to mean Hill's inequality \eqref{eq:Hillsinequality} acting on the Kirchhoff stress $\tau$. In this sense it is used in the FEM-software packages Ansys\textsuperscript{\texttrademark} and Abaqus\textsuperscript{\texttrademark} where, however, no explicit distinction between compressibility and incompressibility is made in Ansys\textsuperscript{\texttrademark}.}
and see Figure \ref{fig:infDruckerstability}, \ref{fig:infDruckerstabilitybox} can thus be expressed as
\begin{align}
	 \frac{1}{2} \frac{\dif^2}{\dif t} \mathcal{E} (t) = \frac{1}{2}\int_{V_0} \langle \dot{\sigma}^{\lin}(t), \dot{\varepsilon}(t) \rangle \dif V_0 > 0
\end{align}
and is sufficient for having a stable equilibrium upon localization. It expresses nothing else than $\C \in \Sym^{++}_4(6)$ for the constitutive law $\sigma = \C . \varepsilon$. Hence, ``positive second order internal work`` is guaranteed by the rate condition $\langle \dot{\sigma} , \dot{\varepsilon} \rangle > 0$, formally similar to the (CSP) requirement \eqref{eq:CSP}.
\subsection{Nonlinear elasticity}
In nonlinear elasticity, the correspondence of positive second order internal work with (CSP) is lost in general. Following, we will show this in more detail:
let us write the internal stored energy in nonlinear elasticity as
\begin{equation}
	\mathcal{E}(t) = \int_{V_0} \WW(F(t)) \dif V_0 \, .
\end{equation}
Then we calculate (recall $\sigma = \frac{1}{J} \, S_1 \, F^T$ and $D = \sym L = \sym (\dot F \, F^{-1}), L= \dot{F} F^{-1}=\DD_{\xi} v(\xi , t)$)
\begin{equation}
	\label{eqPenergy}
	\begin{alignedat}{2}
		\frac{\dif}{\dif t} \mathcal{E}(t) &= \int_{V_0} \langle \DD_F \WW(F(t)), \dot F(t) \rangle \dif V_0 = \int_{V_0} \langle S_1(t), \dot F(t) \rangle \dif V_0 = \int_{V_0}\langle S_1(t) , \dot{F}(t) \, F^{-1} F \rangle \dif V_0\\
		&= \int_{V_0} \langle S_1 \, F^T, \dot F \, F^{-1} \rangle \dif V_0
		= \int_{V_0} \langle \sigma, L \rangle \, \underbrace{J \, \dif V_0}_{= \, \dif V_t} \overset{\sigma \in \Sym(3)}{=} \int_{V_t} \langle \sigma, D \rangle \dif V_t = \int_{V_0} \langle \tau , D \rangle \dif V_0 = \mathcal{P}_{\textnormal{int}}
	\end{alignedat}
\end{equation}
and  with $\frac{\DD}{\DD t}$ denoting the material derivative we obtain
\begin{align}
		\label{eq_secondorderwork}
		\frac{\dif^2}{\dif t^2} \mathcal{E} (t) &= \frac{\dif}{\dif t} \int_{V_0} \langle \sigma (t) , D(t) \rangle \cdot J(t) \dif V_0 \quad \text{(calculate time derivatives only with respect} \nonumber \\[-0.8em]
		&\hspace*{14.5em} \text{ to the fixed referential domain)}  \nonumber \\
		&= \int_{V_0} \frac{\DD}{\DD t}\big( \langle \sigma (t) , D(t) \rangle \cdot J(t)\big) \dif V_0 \nonumber \\
		&= \int_{V_0} \langle \frac{\DD}{\DD t} [\sigma] , D(t) \rangle \cdot J(t) + \langle \sigma , \frac{\DD}{\DD t}[D (t)] \rangle \cdot J(t) + \langle \sigma , D \rangle \cdot \frac{\DD}{\DD t} J(t) \dif V_0 \nonumber \\
		&= \int_{V_0} \langle \frac{\DD}{\DD t} [\sigma] , D(t) \rangle \cdot J(t) + \langle \sigma , \frac{\DD}{\DD t}[D (t)] \rangle \cdot J(t) + \langle \sigma , D \rangle \cdot \langle  \Cof F(t) , \dot{F} (t)\rangle \dif V_0 \nonumber \\
		&= \int_{V_0} \langle \frac{\DD}{\DD t} [\sigma] , D(t) \rangle \cdot J(t) + \langle \sigma , \frac{\DD}{\DD t}[D (t)] \rangle \cdot J(t) + \langle \sigma , D \rangle \cdot \det F(t) \langle F^{-T} (t) , \dot{F} (t) \rangle \dif V_0 \nonumber \\
		&= \int_{V_0} \langle \frac{\DD}{\DD t} [\sigma] , D(t) \rangle \cdot J(t) + \langle \sigma , \frac{\DD}{\DD t}[D (t)] \rangle \cdot J(t) + \langle \sigma , D \rangle \cdot J \langle \id , L \rangle \dif V_0 \nonumber \\
		&= \int_{V_0} \langle \frac{\DD}{\DD t} [\sigma] , D(t) \rangle \cdot J(t) + \langle \sigma , \frac{\DD}{\DD t}[D (t)] \rangle \cdot J(t) + \langle \sigma , D \rangle \cdot J \langle \id , D \rangle \dif V_0\\
		&= \int_{V_0} \big[\langle \frac{\DD}{\DD t} [\sigma] , D(t) \rangle + \langle \sigma , \frac{\DD}{\DD t}[D (t)] \rangle + \langle \sigma , D \rangle \tr (D) \big] \, \underbrace{J (t)\dif V_0}_{\dif V_t} \nonumber \\
		&= \int_{V_t} \big[\langle \frac{\DD}{\DD t} [\sigma] , D(t) \rangle + \langle \sigma , \frac{\DD}{\DD t}[D (t)] \rangle + \langle \sigma , D \rangle \tr (D) \big] \, \dif V_t \nonumber \\
		&\overset{\mathclap{\sigma \in \Sym(3)}}{=} \quad \int_{V_t} \big[\langle \frac{\DD}{\DD t} [\sigma] , D(t) \rangle + \langle \sigma , \frac{\DD}{\DD t}[\DD_{\xi} v(\xi,t)] \rangle + \langle \sigma , D \rangle \tr (D) \big] \, \dif V_t \nonumber \\
		&\overset{\mathclap{\substack{\text{material time} \\ \text{derivative}\vspace*{0.3em}}}}{=} \hspace*{1.5em} \int_{V_t} { \big[\langle \frac{\DD}{\DD t} [\sigma] , D (t) \rangle + \langle \sigma (t) , \DD_{\xi}^2 v(\xi , t) . v + \DD_{\xi} v_{,t}(\xi , t) \cdot 1 \rangle  + \langle \sigma , D \rangle \tr (D) \big]} \dif V_t \nonumber \\
		&= \int_{V_t} \big[ \underbrace{\langle \frac{\DD}{\DD t}[\sigma] , D(t) \rangle}_{\textnormal{not objective!}}
		+ \langle \sigma , \DD_{\xi}^2 v(\xi, t) . v \rangle
		- \langle \underbrace{\textnormal{Div}_{\xi} \, \sigma(t)}_{\substack{= \; 0 \; \textnormal{in spatial} \nonumber \\ \textnormal{equilibrium} \\ \textnormal{(e.g. in homo-} \\ \textnormal{geneous tests)}}}, v_{,t}(\xi,t) \rangle
		+ \underbrace{\div (\sigma^T . v_{,t})}_{\substack{\textnormal{Gauß} \\ \int_{\partial V_t} \langle \sigma^T . v_{,t} ,\overset{\rightarrow}{n} \rangle \dif S_t \\ = \int_{\partial V_t} \langle \sigma . \overset{\rightarrow}{n} , v_{,t} \rangle \dif S_t}} + \langle \sigma , D \rangle \tr (D) \big] \dif V_t \nonumber 
\end{align}
where we used that
\begin{equation}
	\begin{alignedat}{2}
		\div (\sigma^T . v_{,t}) &= \langle \div \sigma , v_{,t} \rangle + \langle \sigma , \DD_{\xi} v_{,t} \rangle \quad
		\implies \quad \langle \sigma, \DD_{\xi} v_{,t} \rangle &= \div (\sigma^T . v_{,t}) - \langle \div \sigma , v_{,t} \rangle \, .
	\end{alignedat}
\end{equation}
It follows
\begin{equation}
	\label{eq:secondorderspatial}
	\begin{alignedat}{2}
		\frac{\dif^2}{\dif t^2} \mathcal{E}(t) &= \int_{V_t} \langle \frac{\DD}{\DD t} [\sigma] , D(t) \rangle \dif V_t 
		+ \int_{V_t} \langle \sigma , \DD_{\xi}^2 v(\xi, t) . v \rangle \dif V_t\\
		&\hspace*{2em}+ \int_{V_t} \langle \sigma , D (t) \rangle \tr (D(t)) \dif V_t + \int_{\partial V_t} \langle \sigma . \overset{\rightarrow}{n} , v_{,t} \rangle \dif S_t
		\equalscolon \frac{\dif^2}{\dif t^2} \mathcal{E}^{\text{spatial}}(t) \, .
	\end{alignedat}
\end{equation}
Thus, the positive ``second order internal work'' criterion $\frac{\dif^2}{\dif t^2} \mathcal{E}^{\text{spatial}}(t) > 0$ in equilibrium for finite strain and motions with $v_{,t}=0$ or $\sigma . \overset{\rightarrow}{n}=0$ at the spatial boundary can be written as the requirement
\begin{align}
	\label{eq:zweiteAblEspatial}
	\frac{\dif^2}{\dif t^2} \mathcal{E}^{\text{spatial}}(t) = \int_{V_t} \langle \frac{\DD}{\DD t} [\sigma(t)], D(t) \rangle + \langle \sigma , D(t) \rangle \tr (D(t)) \dif V_t + \int_{V_t} \langle \sigma , \DD_{\xi}^2 v(\xi, t) . v \rangle \dif V_t > 0 \, .
\end{align}
If we assume in addition $\DD_{\xi}^2 v(\xi, t) = 0$ (i.e. $v(\xi, t)$ is affine in $\xi$, or $L=\DD_{\xi}v(\xi , t)$ is independent of $\xi$), then the second integral cancels and we are left with
\begin{equation}
	\frac{\dif^2}{\dif t^2} \mathcal{E}_{\text{affin}}^{\text{spatial}} (t) = \int_{V_t} \langle \frac{\DD}{\DD t} [\sigma] , D(t) \rangle + \langle \sigma , D(t) \rangle \tr (D(t)) \dif V_t = \int_{V_t} \langle \underbrace{\frac{\DD }{\DD t} [\sigma] + \sigma \tr (D(t))}_{\text{see}\footnotemark} , D(t) \rangle \dif V_t \, .
\end{equation}
\footnotetext{Note the similarity to the non-corotational Biezeno-Hencky stress rate \cite{biezeno1928}, sometimes also called Hill-rate (cf.~\cite{korobeynikov2023}))
\begin{equation}
	\begin{alignedat}{2}
		\frac{\DD^{\text{Hencky}}}{\DD t}[\sigma] &\colonequals &&\; \frac{\DD}{\DD t}[\sigma] + \sigma \, W - W \, \sigma + \sigma \, \tr(D) = \frac{\DD^{\ZJ}}{\DD t}[\sigma] + \sigma \tr(D) \notag \, .\\
	\end{alignedat}
\end{equation}
}
In the incompressible case, $\tr (D)=0$, the latter ``nearly'' coincides with the expression ($\frac{\DD^{\circ}}{\DD t}$ is any corotational rate) 
\begin{align}
	\int_{V_t} \underbrace{\langle \frac{\DD^{\circ}}{\DD t}[\sigma(t)], D(t) \rangle}_{\textnormal{objective!}} \dif V_t > 0 \, .
\end{align}
In nonlinear elasticity, therefore, the local corotational stability requirement
\begin{align}
	\langle \frac{\DD^{\circ}}{\DD t}[\sigma(t)], D(t) \rangle > 0 \qquad \forall \, D \in \Sym(3) \setminus \{0\}
\end{align}
must be clearly distinguished from $\frac{1}{2} \frac{\dif^2}{\dif t^2} \mathcal{E}(t) > 0$, the positive second order internal work, contrary to the geometrically linear case.\\

However, it is quite illuminating and useful to observe that any corotational rate
\begin{equation}
	\frac{\DD^{\circ}}{\DD t}[\sigma (t)] = \frac{\DD}{\DD t} [\sigma] + \sigma \, \Omega^{\circ} - \Omega^{\circ} \sigma \, , \quad \Omega^{\circ} \in \mathfrak{so}(3)
\end{equation}
reduces to $\frac{\DD}{\DD t} [\sigma]$ in situations where the spatially homogeneous $F(t)$ remains diagonal (no rotation effects in the rate), i.e. for constant principal axes (cf.~\cite{poscor2024}). Since \cite{CSP2024,leblond1992constitutive,Leblond2024}
\begin{equation}
	\langle \frac{\DD^{\circ}}{\DD t} [\sigma] , D \rangle > 0 \quad \iff \quad \text{TSTS-M}^{++}
\end{equation}
we infer presently that
\begin{equation}
	\text{TSTS-M}^{++} \quad \implies \quad \langle \frac{\DD}{\DD t} [\sigma] , D \rangle > 0 \quad \quad {\iff \quad \quad \left\langle {{\partial _t}[\sigma],D} \right  \rangle   > 0}
\end{equation}
in homogeneous tests: $\dot{F} \, F^{-1} = L(t)=\DD_{\xi}v(\xi,t)$ is independent of $\xi$ and diagonal.
We will use this implication in section \ref{sec:One_constitutive_condition}.\\

We can repeat the above reasoning starting again from \eqref{eqPenergy} but remaining entirely in the referential domain. This yields
\begin{equation}
	\begin{alignedat}{2}
		\label{eq:secondorderreferential}
		\frac{\dif^2}{\dif t^2} \mathcal{E} (t) &= \int_{V_0} \langle \DD_F S_1 (F) . \dot{F} , \dot{F} \rangle + \langle S_1(F) , F_{,tt} \rangle \dif V_0
		= \int_{V_0} \langle \DD_F^2 \WW (F) . \dot{F} , \dot{F} \rangle + \langle \DD_F \WW(F) , F_{,tt} \rangle  \dif V_0\\
		&= \int_{V_0} \langle \DD_F^2 \WW (F) . \dot{F} , \dot{F} \rangle + \langle \DD_F \WW(F) , \DD \varphi_{,tt} \rangle  \dif V_0\\
		&= \int_{V_0} \langle \DD_F^2 \WW (F) . \dot{F} , \dot{F} \rangle - \langle \underbrace{\div \, \DD_F \WW(F)}_{\substack{= \; 0 \; \textnormal{in referential} \\ \textnormal{equilibrium}}} , \varphi_{,tt} \rangle  \dif V_0 + \int_{\partial V_0} \langle S_1 . n , \varphi_{,tt} (t)  \rangle \dif S_0\\
		&= \int_{V_0} \DD_F^2 \WW (F(t)) . (\dot{F}(t) , \dot{F}(t)) \dif V_0 + \int_{\partial V_0} \langle S_1 . n , \varphi_{,tt} (t)  \rangle \dif S_0 \equalscolon \frac{\dif^2}{\dif t^2} \mathcal{E}^{\text{ref}} (t) \, .
	\end{alignedat}
\end{equation}
For the ``acceleration'' in the material versus the spatial configuration it holds
\begin{equation}
	\label{eq:phi_tt}
	\varphi_{,tt}(x,t) = \frac{\DD}{\DD t} [v(\xi,t)] = \DD_{\xi}v(\xi ,t). v + v_{,t} \, .
\end{equation}
	While clearly by construction
	\begin{equation}
		\begin{alignedat}{2}
			\frac{\dif^2}{\dif t^2} \mathcal{E}^{\textnormal{ref}} (t) &= \frac{\dif^2}{\dif t^2} \mathcal{E}^{\textnormal{spatial}} (t) \, ,\\
		\end{alignedat}
	\end{equation}
	it is not possible to obtain the equivalence
	\begin{equation}
		\int_{V_0}{\DD_F^2 \WW (F(t)) . (\dot{F}(t), \dot{F}(t))} \dif V_0 = \int_{V_t}{ \langle \frac{\DD}{\DD t} [\sigma] + \sigma \tr (D(t)) , D(t) \rangle} \dif V_t + \int_{V_t} \langle \sigma , \DD_{\xi}^2 v(\xi, t) . v \rangle \dif V_t\, ,
	\end{equation}
	since $\varphi_{,tt}=0$ at $\partial V_0 \notiff v_{,t}=0$ at $\partial V_t$ as \eqref{eq:phi_tt} clearly shows.\\

Summarizing, we observe the concordance in quasistatic loading
\begin{align*}
	\textnormal{linear elasticity} & &&  & & \qquad \textnormal{nonlinear elasticity} \\
	\int_{V_0} \frac{1}{2} \langle \sigma, \varepsilon \rangle \dif V_0 \quad & &&  & & \rightsquigarrow \quad \int_{V_0} \WW(F(t)) \dif V_0 \qquad \textnormal{``energy/work'' (objective)} \\
	\int_{V_0} \langle \sigma, \dot{\varepsilon} \rangle \dif V_0 \quad & &&  & &\rightsquigarrow \quad \int_{V_t} \langle \sigma, D \rangle \dif V_t = \int_{V_0} \langle S_1 , \dot{F} \rangle \dif V_0= \mathcal{P}_{\textnormal{int}} \qquad \textnormal{``internal power'', ``rate of work''} \\[-0.8em]
	& && && \hspace*{25em}\textnormal{(objective)}\\
	\int_{V_0} \langle \dot{\sigma}, \dot{\varepsilon} \rangle \dif V_0 \quad & &&  & &\rightsquigarrow \quad \int_{V_t} \langle \frac{\DD}{\DD t}[\sigma] + \sigma \tr(D(t)) , D \rangle + \langle \sigma , \DD_{\xi}^2 v(\xi, t) . v \rangle \dif V_t + \int_{\partial V_t} \langle \sigma . \overset{\rightarrow}{n} , v_{,t} \rangle \dif S_t\\
	&  && && \hspace*{2em}= \underbrace{\int_{V_0} \DD_F^2 \WW (F) . (\dot{F} , \dot{F}) \dif V_0}_{\substack{\textnormal{``second order internal work''} \\ \textnormal{(not objective)}}} + \int_{\partial V_0} \langle S_1 . n , \varphi_{,tt} (t)  \rangle \dif S_0 \, .
\end{align*}
\section{Incremental elastic moduli}
There is no universally accepted unique definition of incremental elastic moduli (cf.~\cite{mihai2017,scott2006}), the only obvious requirement is that they should reduce to their linear elastic counterpart at zero loads. Since all stress tensors coincide to first order in the stress free reference configuration $\id$, different moduli connected to different stress tensors may be defined. Here, we will restrict attention to the Cauchy stress.
For example in Scott \cite{scott2006} we have the definition for the incremental Young's modulus
\begin{equation}
	E_{\text{Scott}}^{\text{incr}} (\lambda_1) \colonequals \lambda_1 \, \DD_{\lambda_1}\widetilde{\sigma}(\lambda_1)
\end{equation}
in uniaxial tension evaluated for the Cauchy stress $\widehat{\sigma}(\lambda_1)$.
We observe by setting $\widehat{\sigma} (\log \lambda_1) \colonequals \widetilde{\sigma} (\lambda_1)$ that
\begin{align}
		\DD_{\lambda_1} \widetilde{\sigma} (\lambda_1) &= \DD_{\lambda_1} [\widehat{\sigma} (\log \lambda_1)] = \DD_{\log \lambda_1} \widehat{\sigma} (\log \lambda_1) \cdot \frac{1}{\lambda_1} \nonumber\\
		& \implies \lambda_1 \, \DD_{\lambda_1} \widetilde{\sigma} (\lambda_1) = \DD_{\log \lambda_1} \widehat{\sigma} (\log \lambda_1) \equalscolon E_{\log}^{\text{incr}} (\lambda_1) \overset{!}{=} E_{\text{Scott}}^{\text{incr}} (\lambda_1) \, .
		\label{eq:Eincrlog}
\end{align}
In \cite{mihai2017} also $E_{\log}^{\text{incr}}$ is used.
In the following we adopt the more primitive definition
\begin{equation}
	\label{eq:Eincr}
	E^{\text{incr}} (\lambda_1) \colonequals \DD_{\lambda_1} \widetilde{\sigma}(\lambda_1) \, , \qquad \qquad E \colonequals E^{\text{incr}} (\lambda_1) \big\vert_{\lambda_1=1} = \DD_{\lambda_1} \widetilde{\sigma} (1)= \DD_{\log \lambda_1} \widehat{\sigma} (0)
\end{equation}
since we are only interested in the positivity of the moduli and \eqref{eq:Eincr} can be interpreted as the slope of the uniaxial Cauchy tension principal-stretch curve, cf. Figure \ref{fig:incrYoung}.\\
\begin{figure}[h!]
	\begin{center}
		\begin{minipage}[h!]{0.45\linewidth}
			\begin{tikzpicture}
				\begin{axis}[scale=1.5,
					axis x line=middle,axis y line=middle,
					width=.7\linewidth,
					height=0.7\linewidth,
					xlabel={$\lambda_1$},        
					ylabel={$\sigma$},        
					xtick={-2},
					ytick={-6},
					xmin=-1,
					xmax=6,
					ymin=-5
					]
					\addplot[blue, thick, domain=0.4:4.5, samples=\sample]{2*((ln(x))/x) *exp((ln(x)^2)};
					\addplot[red, thick, domain=0.3:2]{2*x-2};
					\node at (axis cs:1.5,1.5){$\red{E}$};
					\addplot[red, thick, domain=3.68:4.3]{2.5*x -5.34};
					\draw[red, dashed] (axis cs:3.68, 3.86) -- (axis cs:4.3,3.86);
					\draw[red, dashed] (axis cs:4.3,3.86) -- (axis cs:4.3,5.4);
					\node at (axis cs: 4.9,4.7){$\red{E^{\text{incr}}}$};
				\end{axis}
			\end{tikzpicture}
			\caption{Incremental Young's modulus $E^{\text{incr}}$ according to the simple definition \eqref{eq:Eincr} and Young's modulus $E$ of the infinitesimal theory.}
			\label{fig:incrYoung}
		\end{minipage}
		\qquad
		\begin{minipage}[h!]{0.45\linewidth}
			\begin{tikzpicture}
				\begin{axis}
					[scale=1.5,
					axis x line=middle,axis y line=middle,
					width=.7\linewidth,
					height=0.7\linewidth,
					xlabel={$\log \lambda_1$},        
					ylabel={$\widehat{\sigma}$},        
					xtick={-5},
					ytick={-11},
					xmin=-4,
					xmax=4.5,
					ymin=-10,
					ymax=10
					]
					\addplot[blue, thick, domain=-2:2, samples=\sample]{2*(x/exp(x)) *exp((x^2)};
					\addplot[red, thick, domain=-0.7:0.7, samples=\sample]{2*x};
					\addplot[red, thick, domain=1.3:1.7, samples=\sample]{9.0982*x-7.9875};
					\node at (axis cs: 1.2,1.1){$\red{E}$};
					\draw[red,dashed] (axis cs: 1.3,3.8401) -- (axis cs:1.7,3.8401);
					\draw[red,dashed] (axis cs:1.7,3.8401) -- (axis cs: 1.7,7.47944);
					\node at (axis cs: 2.7,5.15){$\red{E_{\log}^{\text{incr}}}$};
				\end{axis}
			\end{tikzpicture}
			\caption{Incremental Young's modulus in the logarithmic representation according to \eqref{eq:Eincrlog}. The physical content of both $E^{\text{incr}}$ and $E^{\text{incr}}_{\log}$ coincides, only their numerical values differ.}
		\end{minipage}
	\end{center}
\end{figure}

It is clear that
\begin{equation}
	E^{\text{incr}} > 0 \quad \iff \quad E_{\log}^{\text{incr}} > 0 \, ,
\end{equation}
a property that also pertains to the other considered incremental moduli.
\section{One constitutive condition to rule them all: positive incremental moduli for (CSP)}
\label{sec:One_constitutive_condition}
We are now showing that (CSP) in conjunction with a diagonal, homogeneous deformation family $t \mapsto F(t)$ leads to positive incremental moduli in uniaxial tension, equibiaxial extension, planar tension and hydrostatic tension (cf. \cite{Young2024}).
In this respect, let us gather the necessary relations and assumptions:
\begin{equation}
	\begin{alignedat}{4}
		\dot{F} &(t) \, F^{-1} (t)\quad \text{is diagonal} , \ L = \DD_{\xi} v(\xi, t) =\text{const.} \ \text{in} \ \xi \, ,\\
		F(t) &= \diag (\lambda_1 (t) , \lambda_2 (t) , \lambda_3 (t)), \quad
		\dot{F} (t) = \diag (\dot{\lambda_1}(t), \dot{\lambda_2}(t), \dot{\lambda_3}(t)) \, ,\\
		\dot{F} (t) \, F^{-1} (t) &= \diag \big(\frac{\dot{\lambda_1} (t)}{\lambda_1(t)} , \frac{\dot{\lambda_2}(t)}{\lambda_2(t)}, \frac{\dot{\lambda_3}(t)}{\lambda_3(t)}\big) = D (t) \in \Sym(3) \, ,
	\end{alignedat}
\end{equation}
\begin{equation}
	\begin{alignedat}{2}
		&\sigma (V) = \sigma (F) = \diag (\sigma_1(\lambda_1(t), \lambda_2(t), \lambda_3(t)), \sigma_2 (\lambda_1(t), \lambda_2(t), \lambda_3(t)) , \sigma_3 (\lambda_1(t), \lambda_2(t), \lambda_3(t))) \, ,\\
		&\sigma \ \text{is independent of} \ \xi, \ \text{because} \ F(t) \ \text{is homogeneous,}
	\end{alignedat}
\end{equation}
\begin{align}
		\text{(CSP)} \quad &\iff \quad \text{TSTS-M}^{++} \quad \iff 0 < \langle \frac{\DD^{\circ}}{\DD t} [\sigma] , D \rangle \overset{F \, = \, \diag}{=} \langle \frac{\DD}{\DD t} [\sigma] , D \rangle \, , \nonumber\\
		0 < \langle \frac{\DD}{\DD t} [\sigma] , D \rangle &= \langle \frac{\DD}{\DD t}[\diag (\sigma_1 (\lambda_1(t) , \lambda_2(t), \lambda_3(t)), \sigma_2(\dots), \sigma_3(\dots))] , \diag \big( \frac{\dot{\lambda_1}(t)}{\lambda_1(t)} , \frac{\dot{\lambda_2}(t)}{\lambda_2(t)} , \frac{\dot{\lambda_3}(t)}{\lambda_3(t)}\big) \rangle  \\
		&\overset{\mathclap{\substack{\DD_{\xi} \sigma = 0: \\ \text{material time} \\ \text{derivative} \vspace*{0.3em}}}}{\underset{\mathclap{\frac{\DD}{\DD t} \rightarrow \, \partial_t}}{=}} \quad \big\langle \begin{pmatrix}
			\partial_t [\sigma_1 (\lambda_1(t), \lambda_2(t), \lambda_3(t))]\\
			\partial_t [\sigma_2 (\lambda_1(t), \lambda_2(t), \lambda_3(t))]\\
			\partial_t [\sigma_3 (\lambda_1(t), \lambda_2(t), \lambda_3(t))]
		\end{pmatrix} ,
		\begin{pmatrix}
			\frac{\dot{\lambda_1} (t)}{\lambda_1(t)}\\
			\frac{\dot{\lambda_2} (t)}{\lambda_2(t)}\\
			\frac{\dot{\lambda_3} (t)}{\lambda_3(t)}\\
		\end{pmatrix}
	\big\rangle_{\R^3}
	=  \sum_{i=1}^{3} {\partial_t [\sigma_i(\lambda_1(t), \lambda_2(t), \lambda_3(t))] \cdot \frac{\dot{\lambda_i}(t)}{\lambda_i(t)}} \, \nonumber.
\end{align}
\subsection{Positive incremental moduli for compressible response}
\begin{itemize}
	\item \textbf{compressible uniaxial tension}, i.e. $\sigma_2=\sigma_3=0, \lambda_1$ free, $\lambda_2=\lambda_3=\lambda_2(\lambda_1(t))$ and we have
	\begin{equation}
		\begin{alignedat}{2}
			\label{eq:compuniaxten}
			0 < \langle \frac{\DD}{\DD t} [\sigma] , D \rangle &= \sum_{i=1}^{3} {\partial_t [\sigma_i(\lambda_1(t), \lambda_2(t), \lambda_3(t))] \cdot \frac{\dot{\lambda_i}(t)}{\lambda_i(t)}}\\
			&= \partial_t [\sigma_1(\lambda_1(t), \lambda_2 (\lambda_1(t)), \lambda_2(\lambda_1(t)))] \cdot \frac{\dot{\lambda_1}(t)}{\lambda_1(t)} + 0 + 0\\
			&= \partial_t [\underbrace{\sigma_1(\lambda_1(t), \lambda_2 (\lambda_1(t)), \lambda_2(\lambda_1(t)))}_{\equalscolon \widetilde{\sigma}(\lambda_1(t))}] \cdot \frac{\dot{\lambda_1}(t)}{\lambda_1(t)}\\
			&= \DD_{\lambda_1} \widetilde{\sigma} (\lambda_1(t)) \cdot \dot{\lambda_1}(t) \cdot \frac{\dot{\lambda_1}(t)}{\lambda_1(t)}\\
			&= \DD_{\lambda_1} \widetilde{\sigma} (\lambda_1(t)) \cdot \frac{\vert\dot{\lambda_1}(t)\vert^2}{\lambda_1(t)} > 0  \quad 
			\implies \DD_{\lambda_1}\widetilde{\sigma} (\lambda_1(t)) > 0 \, .
		\end{alignedat}
	\end{equation}
	Thus $\lambda_1 \mapsto \widetilde{\sigma}(\lambda_1)$ is monotone and the incremental Young's modulus $E^{\text{incr}}$ is positive, \\ $E^{\text{incr}}\colonequals \DD_{\lambda_1} \widetilde{\sigma}(\lambda_1) > 0$.
	\item \textbf{compressible equibiaxial extension}, i.e. $\lambda_1=\lambda_2, \sigma_1 =\sigma_2$ and $\sigma_3=0$
	\begin{align}
			0 < \langle \frac{\DD}{\DD t} [\sigma] , D \rangle &= \sum_{i=1}^{3} {\partial_t [\sigma_i(\lambda_1(t), \lambda_2(t), \lambda_3(t))] \cdot \frac{\dot{\lambda_i}(t)}{\lambda_i(t)}} \nonumber \\
			&= \partial_t [\sigma_1 (\lambda_1(t) , \lambda_1 (t), \lambda_3(\lambda_1(t))] \cdot \frac{\dot{\lambda_1} (t)}{\lambda_1(t)} + \partial_t [\sigma_1 (\lambda_1(t) , \lambda_1 (t), \lambda_3(\lambda_1(t))] \cdot \frac{\dot{\lambda_2} (t)}{\lambda_2(t)} + 0  \nonumber \\
			&= \partial_t [\sigma_1 (\lambda_1(t) , \lambda_1 (t), \lambda_3(\lambda_1(t))] \cdot \frac{\dot{\lambda_1} (t)}{\lambda_1(t)} + \partial_t [\sigma_1 (\lambda_1(t) , \lambda_1 (t), \lambda_3(\lambda_1(t))] \cdot \frac{\dot{\lambda_1} (t)}{\lambda_1(t)} \nonumber \\
			&= 2 \cdot \partial_t [\underbrace{\sigma_1 (\lambda_1(t) , \lambda_1 (t), \lambda_3(\lambda_1(t))}_{\equalscolon \widetilde{\sigma}(\lambda_1(t))}] \cdot \frac{\dot{\lambda_1} (t)}{\lambda_1(t)}\\
			&= 2 \cdot \DD_{\lambda_1} \widetilde{\sigma} (\lambda_1(t)) \cdot \dot{\lambda_1}(t) \cdot \frac{\dot{\lambda_1}(t)}{\lambda_1(t)} \nonumber \\
			&= 2 \cdot \DD_{\lambda_1} \widetilde{\sigma} (\lambda_1(t)) \cdot \frac{\vert \dot{\lambda_1}(t)\vert^2}{\lambda_1(t)} > 0 \quad
			\implies \DD_{\lambda_1}\widetilde{\sigma} (\lambda_1(t)) > 0 \, . \nonumber 
	\end{align}
	Thus $\lambda_1 \mapsto \widetilde{\sigma}(\lambda_1)$ is monotone and the incremental equibiaxial modulus $A^{\text{incr}}$ is positive, \\ $A^{\text{incr}} \colonequals \frac{1}{2} \DD_{\lambda_1} \widetilde{\sigma}(\lambda_1)> 0$.
	\item \textbf{compressible planar tension}, i.e. $\lambda_3=1$ and $\sigma_2=0$
	\begin{align}
			0 < \langle \frac{\DD}{\DD t} [\sigma] , D \rangle &= \sum_{i=1}^{3} {\partial_t [\sigma_i(\lambda_1(t), \lambda_2(t), \lambda_3(t))] \cdot \frac{\dot{\lambda_i}(t)}{\lambda_i(t)}}  \nonumber \\
			&= \partial_t [\sigma_1(\lambda_1(t) , \lambda_2(\lambda_1 (t)) , 1)] \cdot \frac{\dot{\lambda_1} (t)}{\lambda_1 (t)} + 0 + \partial_t [\sigma_3 (\lambda_1(t) , \lambda_2(\lambda_1(t)), 1)] \cdot \frac{\dot{\lambda_3}}{\lambda_3} \nonumber \\
			&= \partial_t [\sigma_1(\lambda_1(t) , \lambda_2(\lambda_1 (t)) , 1)] \cdot \frac{\dot{\lambda_1} (t)}{\lambda_1 (t)} + \partial_t [\sigma_3 (\lambda_1(t) , \lambda_2(\lambda_1(t)), 1)] \cdot \underbrace{\frac{\dot{\lambda_3}}{\lambda_3}}_{= 0} \nonumber \\
			&= \partial_t [\underbrace{\sigma_1(\lambda_1(t) , \lambda_2(\lambda_1 (t)) , 1)}_{\equalscolon \widetilde{\sigma}(\lambda_1(t))}] \cdot \frac{\dot{\lambda_1} (t)}{\lambda_1 (t)}\\
			&= \DD_{\lambda_1} \widetilde{\sigma} (\lambda_1(t)) \cdot \dot{\lambda_1}(t) \cdot \frac{\dot{\lambda_1}(t)}{\lambda_1(t)} \nonumber \\
			&= \DD_{\lambda_1} \widetilde{\sigma} (\lambda_1(t)) \cdot \frac{\vert \dot{\lambda_1}(t)\vert^2}{\lambda_1(t)} > 0 \quad
			\implies \DD_{\lambda_1} \widetilde{\sigma} (\lambda_1 (t)) > 0 \, . \nonumber 
	\end{align}
	Thus $\lambda_1 \to \widetilde{\sigma} (\lambda_1)$ is monotone and the incremental planar tension modulus $PT^{\text{incr}}$ is positive, \\ $PT^{\text{incr}} \colonequals \DD_{\lambda_1} \widetilde{\sigma} (\lambda_1) > 0$.
	\item \textbf{compressible hydrostatic tension} (equitriaxial tension), i.e. $\lambda_1(t) = \lambda_2(t) = \lambda_3(t), \ \sigma_1 = \sigma_2 = \sigma_3$
	\begin{align}
			0 < \langle \frac{\DD}{\DD t} [\sigma] , D \rangle &= \sum_{i=1}^{3} {\partial_t [\sigma_i(\lambda_1(t), \lambda_2(t), \lambda_3(t))] \cdot \frac{\dot{\lambda_i}(t)}{\lambda_i(t)}} \nonumber \\
			&= \partial_t [\sigma_1(\lambda_1(t) , \lambda_1(t) , \lambda_1(t))] \cdot \frac{\dot{\lambda_1} (t)}{\lambda_1 (t)} + \partial_t [\sigma_1(\lambda_1(t) , \lambda_1(t) , \lambda_1(t))] \cdot \frac{\dot{\lambda_1} (t)}{\lambda_1 (t)}  \nonumber \\
			&\hspace*{1em}+ \partial_t [\sigma_1(\lambda_1(t) , \lambda_1(t) , \lambda_1(t))] \cdot \frac{\dot{\lambda_1} (t)}{\lambda_1 (t)}  \nonumber \\
			&= 3 \cdot \partial_t [\underbrace{\sigma_1(\lambda_1(t) , \lambda_1(t) , \lambda_1(t))}_{\equalscolon \widetilde{\sigma}(\lambda_1(t))}] \cdot \frac{\dot{\lambda_1} (t)}{\lambda_1 (t)}\\
			&= 3 \cdot \DD_{\lambda_1} \widetilde{\sigma} (\lambda_1(t)) \cdot \dot{\lambda_1}(t) \cdot \frac{\dot{\lambda_1}(t)}{\lambda_1(t)} \nonumber \\
			&= 3 \cdot \DD_{\lambda_1} \widetilde{\sigma} (\lambda_1(t)) \cdot \frac{\vert \dot{\lambda_1}(t)\vert^2}{\lambda_1(t)} > 0 \quad
			\implies \DD_{\lambda_1} \widetilde{\sigma} (\lambda_1(t)) > 0 \, . \nonumber 
	\end{align}
	Thus $\lambda_1 \mapsto \widetilde{\sigma} (\lambda_1)$ is monotone and the incremental bulk modulus $\kappa^{\text{incr}}$ is positive, \\ $\kappa^{\text{incr}} \colonequals \frac{1}{3}\DD_{\lambda_1} \widetilde{\sigma} (\lambda_1(t))
	> 0$.
\end{itemize}
\subsection{Positive incremental moduli for incompressible response and Hill's inequality}
\label{subsec:posincrmoduliincomp}
In the incompressible case ($\det F = 1$), the Cauchy stress tensor $\sigma$ can be identified with the Kirchhoff stress tensor $\tau= \det F \cdot \sigma$.
Moreover, incompressibility allows to nicely by-pass the stress conditions at the free surface and the constitutive law. Instead one can immediately arrive at the kinematics for the family $t \mapsto F(t)$. Only in a second step do we need to calculate the indetermined pressure $p$ from the boundary conditions in order to finally obtain the explicit form of the principal stresses. We will see that Hill's inequality (applied to the incompressible response) is equivalent to (CSP) and already implies positive incremental moduli.
\begin{itemize}
	\item \textbf{incompressible uniaxial tension}, i.e. $\lambda_1 \, \lambda_2 \, \lambda_3 = 1, \ \lambda_2 = \lambda_3, \ \lambda_2 = \frac{1}{\sqrt{\lambda_1}}, \quad (\tau_2=\tau_3=0)$
	\begin{equation}
		\begin{alignedat}{2}
			\label{eq:incunitens}
			\underbrace{0 <  \langle \frac{\DD}{\DD t} [\sigma] , D \rangle}_{\text{from (CSP)}} = \underbrace{\langle \frac{\DD}{\DD t} [\tau] , D \rangle}_{\text{Hill's inequality}} &= \sum_{i=1}^{3} {\partial_t [\tau_i(\lambda_1(t), \frac{1}{\sqrt{\lambda_1(t)}}, \frac{1}{\sqrt{\lambda_1(t)}})] \cdot \frac{\dot{\lambda_i}(t)}{\lambda_i(t)}}\\
			&= \partial_t [\tau_1(\lambda_1(t), \frac{1}{\sqrt{\lambda_1(t)}}, \frac{1}{\sqrt{\lambda_1(t)}})] \cdot \frac{\dot{\lambda_1}(t)}{\lambda_1(t)} + 0 + 0\\
			&= \partial_t [\underbrace{\tau_1(\lambda_1(t), \frac{1}{\sqrt{\lambda_1(t)}}, \frac{1}{\sqrt{\lambda_1(t)}})}_{\equalscolon \widetilde{\tau}(\lambda_1(t))}] \cdot \frac{\dot{\lambda_1}(t)}{\lambda_1(t)}\\
			&= \DD_{\lambda_1} \widetilde{\tau} (\lambda_1(t)) \cdot \dot{\lambda_1}(t) \cdot \frac{\dot{\lambda_1}(t)}{\lambda_1(t)}\\
			&= \DD_{\lambda_1} \widetilde{\tau} (\lambda_1(t)) \cdot \frac{\vert\dot{\lambda_1}(t)\vert^2}{\lambda_1(t)} > 0  \quad 
			\implies \DD_{\lambda_1}\widetilde{\tau} (\lambda_1(t)) > 0 \, .
		\end{alignedat}
	\end{equation}
	Thus $\lambda_1 \mapsto \widetilde{\tau}(\lambda_1)$ is monotone and the incremental Young's modulus $E^{\text{incr}}$ for incompressible response is positive, $E^{\text{incr}}\colonequals \DD_{\lambda_1} \widetilde{\tau}(\lambda_1) > 0$.
	\item \textbf{incompressible equibiaxial extension}, i.e. $\lambda_1=\lambda_2$ and $\lambda_3 = \frac{1}{\lambda_2^2} \ (\tau_3=0)$
	\begin{align}
			0 < \langle \frac{\DD}{\DD t} [\tau] , D \rangle &= \sum_{i=1}^{3} {\partial_t [\tau_i(\lambda_1(t), \lambda_1(t), \frac{1}{\lambda_1^2(t)})] \cdot \frac{\dot{\lambda_i}(t)}{\lambda_i(t)}} \nonumber \\
			&= \partial_t [\tau_1 (\lambda_1(t), \lambda_1(t), \frac{1}{\lambda_1^2(t)})] \cdot \frac{\dot{\lambda_1} (t)}{\lambda_1(t)}  \nonumber \\
			&\hspace*{1em}+ \partial_t [\tau_2 (\lambda_1(t), \lambda_1(t), \frac{1}{\lambda_1^2(t)})] \cdot \frac{\dot{\lambda_2} (t)}{\lambda_2(t)} + 0 \nonumber \\
			&= \partial_t [\tau_1 (\lambda_1(t), \lambda_1(t), \frac{1}{\lambda_1^2(t)})] \cdot \frac{\dot{\lambda_1} (t)}{\lambda_1(t)} + \partial_t [\tau_1 (\lambda_1(t), \lambda_1(t), \frac{1}{\lambda_1^2(t)})] \cdot \frac{\dot{\lambda_1} (t)}{\lambda_1(t)} \nonumber \\
			&= 2 \cdot \partial_t [\underbrace{\tau_1 (\lambda_1(t), \lambda_1(t), \frac{1}{\lambda_1^2(t)})}_{\equalscolon \widetilde{\tau}(\lambda_1(t))}] \cdot \frac{\dot{\lambda_1} (t)}{\lambda_1(t)}\\
			&= 2 \cdot \DD_{\lambda_1} \widetilde{\tau} (\lambda_1(t)) \cdot \dot{\lambda_1}(t) \cdot \frac{\dot{\lambda_1}(t)}{\lambda_1(t)} \nonumber \\
			&= 2 \cdot \DD_{\lambda_1} \widetilde{\tau} (\lambda_1(t)) \cdot \frac{\vert \dot{\lambda_1}(t)\vert^2}{\lambda_1(t)} > 0 \quad
			\implies \DD_{\lambda_1}\widetilde{\tau} (\lambda_1(t)) > 0 \, . \nonumber 
	\end{align}
	Thus $\lambda_1 \mapsto \widetilde{\tau}(\lambda_1)$ is monotone and the incremental equibiaxial modulus $A^{\text{incr}}$ for incompressible response is positive, $A^{\text{incr}} \colonequals \frac{1}{2} \DD_{\lambda_1} \widetilde{\tau}(\lambda_1)> 0$.
	\item \textbf{incompressible planar tension}, i.e. $\lambda_3=1$ and $\lambda_2 = \frac{1}{\lambda_1}\ (\tau_2=0)$
	\begin{align}
			0 < \langle \frac{\DD}{\DD t} [\tau] , D \rangle &= \sum_{i=1}^{3} {\partial_t [\tau_i(\lambda_1(t), \frac{1}{\lambda_1(t)}, 1)] \cdot \frac{\dot{\lambda_i}(t)}{\lambda_i(t)}} \nonumber \\
			&= \partial_t [\tau_1(\lambda_1(t), \frac{1}{\lambda_1(t)}, 1)] \cdot \frac{\dot{\lambda_1} (t)}{\lambda_1 (t)} + 0 + \partial_t [\tau_3 (\lambda_1(t), \frac{1}{\lambda_1(t)}, 1)] \cdot \frac{\dot{\lambda_3}}{\lambda_3} \nonumber \\
			&= \partial_t [\tau_1(\lambda_1(t), \frac{1}{\lambda_1(t)}, 1)] \cdot \frac{\dot{\lambda_1} (t)}{\lambda_1 (t)} + \partial_t [\tau_3 (\lambda_1(t), \frac{1}{\lambda_1(t)}, 1)] \cdot \underbrace{\frac{\dot{\lambda_3}}{\lambda_3}}_{= 0} \nonumber \\
			&= \partial_t [\underbrace{\tau_1(\lambda_1(t), \frac{1}{\lambda_1(t)}, 1)}_{\equalscolon \widetilde{\tau}(\lambda_1(t))}] \cdot \frac{\dot{\lambda_1} (t)}{\lambda_1 (t)}\\
			&= \DD_{\lambda_1} \widetilde{\tau} (\lambda_1(t)) \cdot \dot{\lambda_1}(t) \cdot \frac{\dot{\lambda_1}(t)}{\lambda_1(t)} \nonumber \\
			&= \DD_{\lambda_1} \widetilde{\tau} (\lambda_1(t)) \cdot \frac{\vert \dot{\lambda_1}(t)\vert^2}{\lambda_1(t)} > 0 \quad
			\implies \DD_{\lambda_1} \widetilde{\tau} (\lambda_1 (t)) > 0 \, .  \nonumber 
	\end{align}
	Thus $\lambda_1 \to \widetilde{\tau} (\lambda_1)$ is monotone and the incremental planar tension modulus $PT^{\text{incr}}$ for incompressible response is positive, $PT^{\text{incr}} \colonequals \DD_{\lambda_1} \widetilde{\tau} (\lambda_1) > 0$.
\end{itemize}
\section{Examples: uniaxial tension}
Let us explain with three examples the content of the above development. We consider the uniaxial tension of a cube of homogeneous and isotropic elastic material. Upon loading the cube deforms homogeneously into a parallelepiped and $F(t)$ remains diagonal, thus rotation effects are absent in the corotational derivative and the principal axes remain constant.
	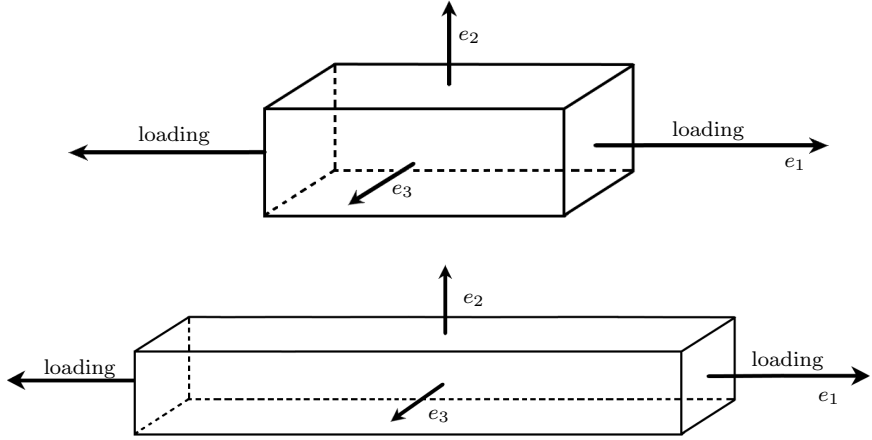
\begin{figure}[h!]
	\begin{center}
		\begin{minipage}[h!]{0.95\linewidth}
			\centering
			\hspace*{40pt}
			\includegraphics[scale=0.35]{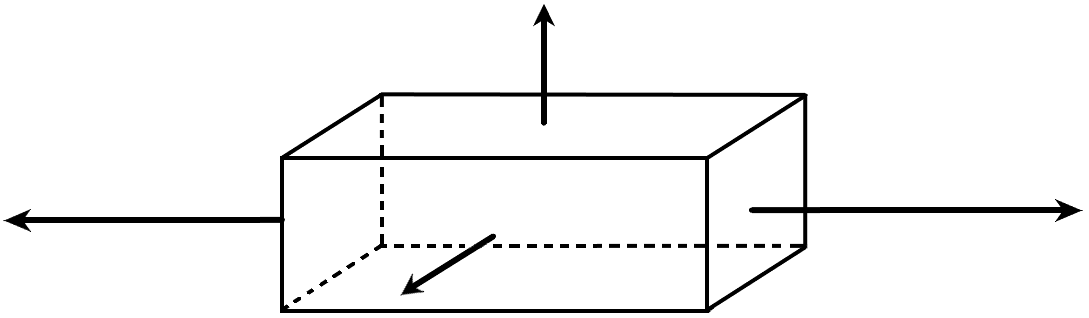}
			\put(-18,20){\footnotesize $e_1$}
			\put(-60,32){\footnotesize loading}
			\put(-260,30){\footnotesize loading}
			\put(-165,10){\footnotesize $e_3$}
			\put(-140,68){\footnotesize $e_2$}
			\vspace*{15pt} \newline
			\centering
			\includegraphics[scale=0.27]{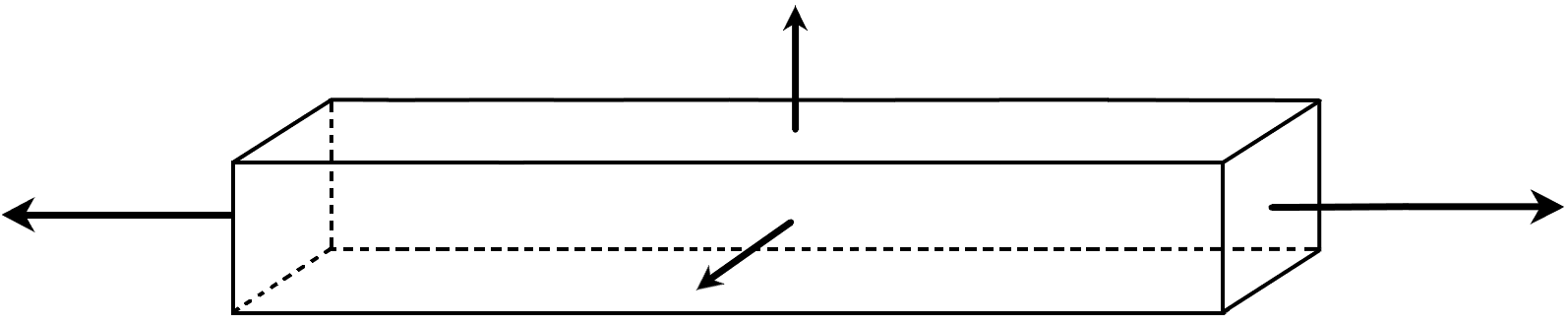}
			\put(-20,15){\footnotesize $e_1$}
			\put(-45,28){\footnotesize loading}
			\put(-310,25){\footnotesize loading}
			\put(-166,8){\footnotesize $e_3$}
			\put(-153,52){\footnotesize $e_2$}
			\caption{Three-dimensional example of a rectangular beam that is pulled in $e_1$-direction. Reference configuration above, deformed configuration below. Uniaxial tension leading to a homogeneous family $F(t)= \diag (\lambda_1 (t) ,\lambda_2 (t) , \lambda_3(t))$.}
			\label{fig1}
		\end{minipage}
	\end{center}
\end{figure}
\subsection{Uniaxial tension: compressible response}
\begin{example}[Exponentiated Hencky energy]
	\label{ex:expHencky_uniaxial}
	The exponentiated Hencky energy satisfies TSTS-M$^{++}$ throughout (cf. \cite{agn_neff2015exponentiatedI,agn_neff2015exponentiatedII}) and is given by
	\begin{align}
		\WW_{\exp} (F)= \frac{\mu}{k} \, \exp\left(k \, \norm{\log V}^2\right) + \frac{\lambda}{2 \, \widehat{k}} \, \exp\left(\widehat{k} \, (\log(\det V))^2\right)
	\end{align}
	with principal Cauchy stresses
	\begin{align}
		\sigma_i(\lambda_1, \lambda_2, \lambda_3) = \frac{1}{\lambda_1 \, \lambda_2 \, \lambda_3} \, \left\{2 \, \mu \, \exp \, (k \, \left(\sum_{j=1}^3 (\log \lambda_j)^2 \right) \, ) \, \log \lambda_i + \lambda \, \exp \left(\widehat{k} \, (\log (\lambda_1 \, \lambda_2 \, \lambda_3))^2 \right) \, \log(\lambda_1 \, \lambda_2 \, \lambda_3)\right\}.
	\end{align}
	For the sake of simplicity, let us assume that $\mu = k = \widehat{k} = 1$ and $\lambda = 2$. Since the lateral sides are free, the equation $\sigma_2(\lambda_1, \lambda_2, \lambda_2) = \sigma_3(\lambda_1, \lambda_2, \lambda_2)= 0$ amounts to
	\begin{align}
		\exp\left( (\log \lambda_1)^2 + 2 \, (\log \lambda_2)^2 - (\log (\lambda_1 \, \lambda_2^2))^2\right) \, \log \lambda_2 + \log (\lambda_1 \, \lambda_2^2) = 0 \cdot
	\end{align} 
		At this point we numerically solve this equation in the form $\lambda_2 = \lambda_2(\lambda_1)$. This solution is then inserted into
		\begin{align}
			\sigma_1(\lambda_1) = \frac{2}{\lambda_1 \, \lambda_2^2} \cdot \left\{ \exp ((\log \lambda_1)^2 + 2 \, (\log \lambda_2)^2) \cdot \log \lambda_1 + \exp \left( (\log(\lambda_1 \, \lambda_2^2))^2\right) \cdot \log(\lambda_1 \, \lambda_2^2) \right\}
		\end{align}
		yielding a plot for $\widetilde {{\sigma}} _1\left( {{\lambda _1}} \right)$. The Biot stress $\widetilde {T}_{\Biot}^1$  is calculated from 
		\begin{eqnarray}
		  \widetilde T_{\text{Biot}}^1({\lambda _1}) = {\lambda _2} \, {\lambda _3} \, {{\widetilde{\sigma}}_1}\left( {{\lambda _1}} \right) = \lambda _2^2\left({\lambda _1} \right){{\widetilde{{\sigma}}}_1} \left({\lambda _1} \right) \cdot
		\end{eqnarray}
Similarly, we evaluate the energy in uniaxial tension and its derivative $\dfrac{\text{d}}{{\text{d}{\lambda _1}}}\widetilde W\left( {{\lambda _1}} \right)$. According to our result in \eqref{eq:compuniaxten} we expect a monotone increasing response $\lambda_1 \mapsto \widetilde{{\sigma}}_1 (\lambda_1)$.
		\begin{figure}[h!]
        \hspace*{\fill}%
            \begin{minipage}[t]{0.45\textwidth}
                \centering
                \vspace{0pt}
                \includegraphics[width=\textwidth]{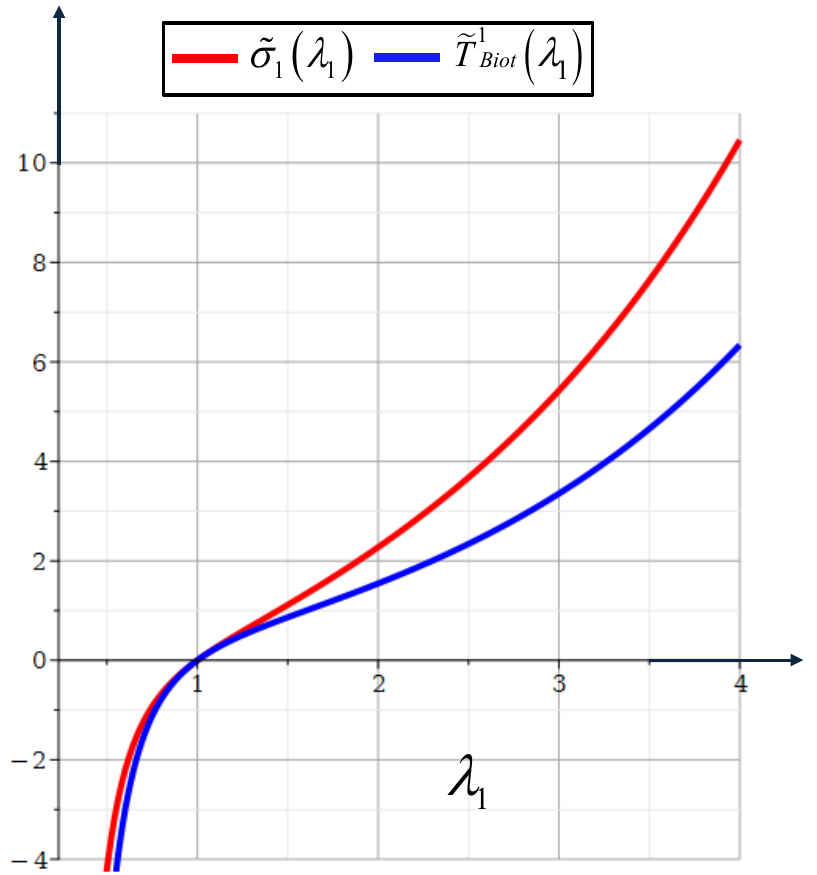}
                \captionof{figure}{Exp-Hencky: monotone tensile Cauchy stress (red) and monotone $T_{\Biot}$ stress (blue).}
                \label{Exple_5_1}
            \end{minipage}%
            \hfill
            \begin{minipage}[t]{0.45\textwidth}
                \centering
                \vspace{0pt}
                \includegraphics[width=\textwidth]{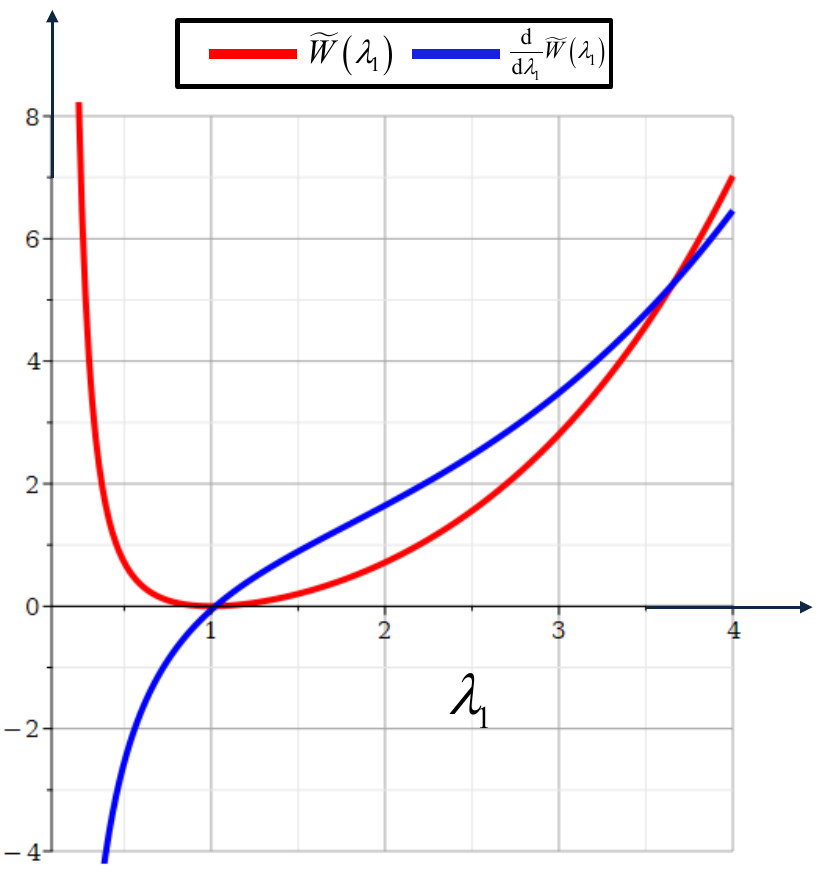}
                \captionof{figure}{Exp-Hencky: convex energy in uniaxial tension and its monotone derivative.}
                \label{Energy_5_1}
            \end{minipage}%
        \hspace*{\fill}
        \end{figure} 
\end{example}
\begin{example}[compressible Neo-Hooke with vol-iso split] 
	\begin{equation}
		\begin{alignedat}{2}
			\label{eq:WNHcomp}
			\WW_{\NH}(F) &= \frac{\mu}{2} \big( \frac{\norm{F}^2}{(\det F)^{\frac{2}{3}}} -3 \big) + \frac{\kappa}{2} (\det F - 1)^2 \, ,
		\end{alignedat}
	\end{equation}
	\begin{equation*}
		\label{eq:sigmaNHcomp}
		\sigma_{\NH} (B) =  \frac{\mu}{(\det B)^{5/6}} \cdot \dev_3 B + \kappa \, (\sqrt{\det B} - 1) \id
		= \frac{\mu}{(\det B)^{5/6}} \cdot (B - \frac{1}{3} \, \tr \, B \cdot \id) + \kappa \, (\sqrt{\det B} - 1) \cdot \id
	\end{equation*}
	This compressible Neo-Hooke model does not satisfy (CSP). The principal Cauchy stresses are given by
	\begin{align}
		\sigma_{\NH}^{i} (\lambda_1, \lambda_2, \lambda_3) = \frac{ \mu}{(\lambda_1^2 \lambda_2^2 \lambda_3^2)^{5/6}} \, (\lambda_i^2 - \frac{1}{3} (\lambda_1^2 + \lambda_2^2 + \lambda_3^2)) + \kappa \, (\lambda_1\lambda_2\lambda_3 - 1) \cdot
	\end{align}
Assuming $\sigma_1 = \sigma_2 = 0$ gives again ${\lambda _2} = {\lambda _3}\left( {{\lambda _1}} \right)$, and we calculate as before the uniaxial stresses  $\widetilde{\sigma}_1$ and $\widetilde{T}_{\Biot}^1$. Our result in \eqref{eq:compuniaxten} is not directly applicable, but the response $\lambda_1 \mapsto \widetilde{{\sigma}}_1 (\lambda_1)$ is nevertheless monotone.
		\begin{figure}[h!]
        \hspace*{\fill}%
            \begin{minipage}[t]{0.45\textwidth}
                \centering
                \vspace{0pt}
                \includegraphics[width=\textwidth]{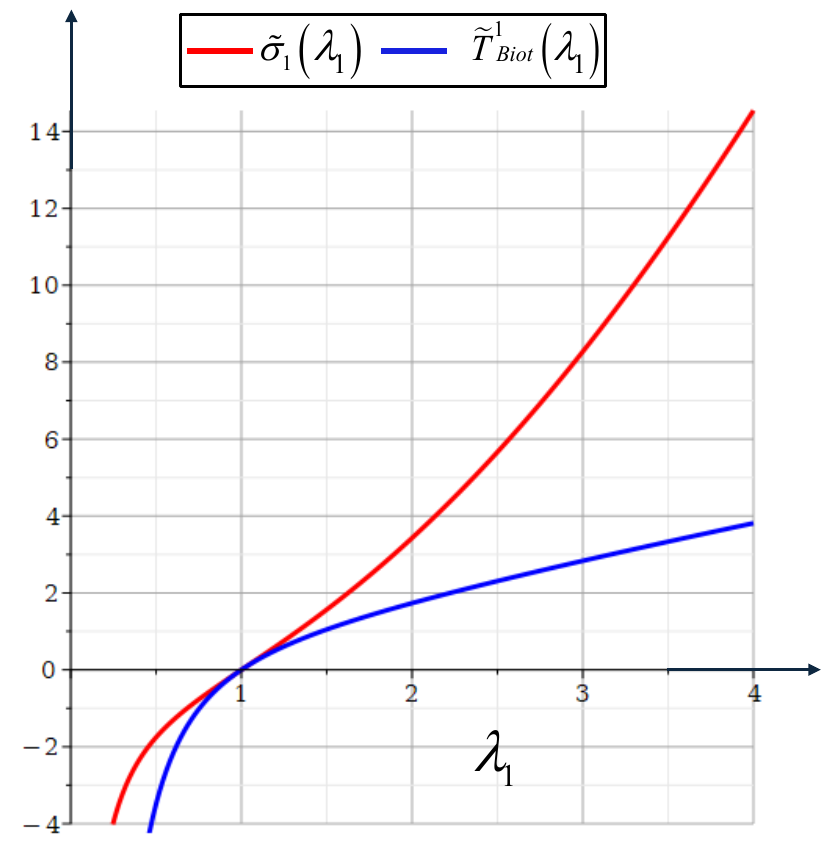}
                \captionof{figure}{Compressible Neo-Hooke: tensile Cauchy stress $\sigma_1$ and $T_{\Biot}$-stress are still both monotone increasing while (CSP) is not satisfied.}
                \label{Exple_5_2}
            \end{minipage}%
            \hfill
            \begin{minipage}[t]{0.45\textwidth}
                \centering
                \vspace{0pt}
                \includegraphics[width=\textwidth]{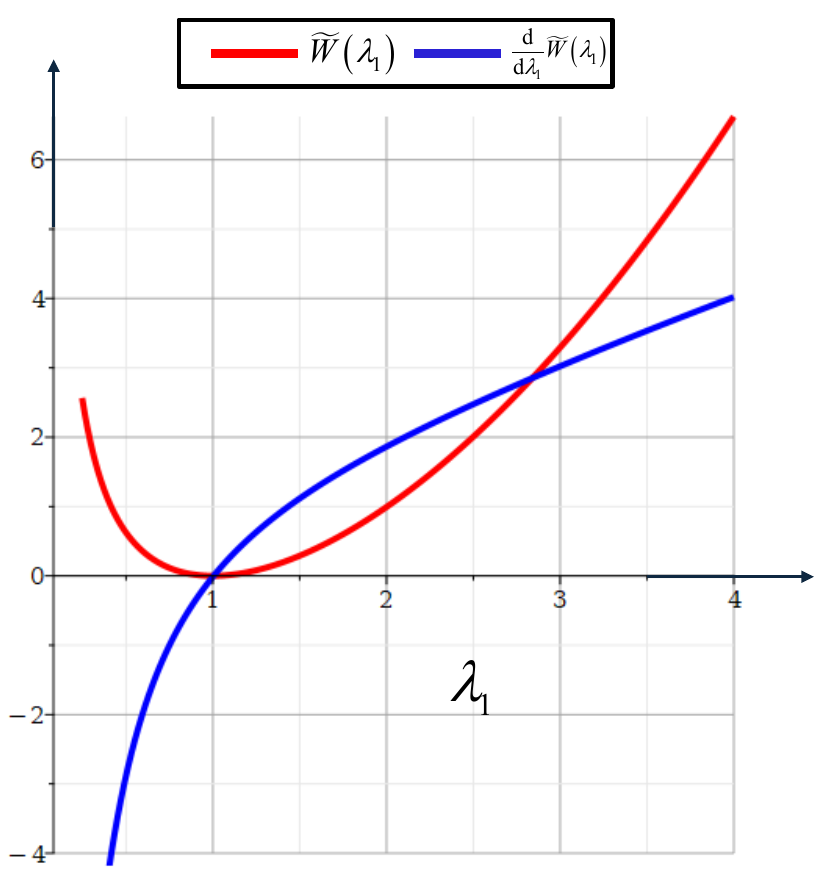}
                \captionof{figure}{Compressible Neo-Hooke: convex uniaxial energy and its monotone derivative.}
                \label{Energy_5_2}
            \end{minipage}%
        \hspace*{\fill}
        \end{figure} 
  \end{example}
\begin{example}[quadratic Hencky energy] \label{ex3.2}
	A very similar example to the first one is given by the non LH-elliptic hyperelastic quadratic Hencky energy (see \cite{Bruhns01, agn_ghiba2015ellipticity, agn_martin2018non, agn_nedjar2018finite, agn_neff2015geometry, agn_neff2014axiomatic, agn_neff2015exponentiatedI, agn_neff2015exponentiatedII} for related literature).\\
	However, the quadratic Hencky energy does not satisfy the (CSP) condition while Hill's inequality holds. We have
	\begin{align}
			\WW(F) &= \mu \, \norm{\log V}^2 + \frac{\lambda}{2} \, \tr^2(\log V),  \quad \quad \tau_{\text{Hencky}} (V) = 2 \, \mu \log V + \lambda \tr (\log V) \cdot \id \, ,\\
			\sigma(V) &= \frac{1}{\det V} \, \left(2 \, \mu \, \log V + \lambda \, \log \det V \, \id\right)
			= \frac{1}{\det V} \, \left( \frac{E}{(1 + \nu)} \, \log V + \frac{E \, \nu}{(1+\nu) \, (1-2 \, \nu)} \, \log \det V \, \id\right) \nonumber
	\end{align}
	with principal stresses
	\begin{align}
		\sigma_i &= \frac{1}{\lambda_1 \, \lambda_2 \, \lambda_3} \, \{2 \, \mu \, \log \lambda_i + \lambda \cdot \log(\lambda_1 \, \lambda_2 \, \lambda_3) \} = \frac{1}{\lambda_1 \, \lambda_2 \, \lambda_3} \, \log \left(\lambda_i^{2 \, \mu} \, (\lambda_1 \, \lambda_2 \, \lambda_3)^{\lambda} \right)\\
		&=\frac{1}{\lambda_1\, \lambda_2 \, \lambda_3} \log \big( \lambda_i^{2 \, \frac{E}{2 (1 + \nu)}} (\lambda_1 \, \lambda_2 \, \lambda_3)^{\frac{E \, \nu}{(1 + \nu) (1-2 \, \nu)}}\big)
		=\frac{1}{\lambda_1\, \lambda_2 \, \lambda_3} \log \big( \lambda_i^{\frac{E}{1 + \nu}} (\lambda_1 \, \lambda_2 \, \lambda_3)^{\frac{E \, \nu}{(1 + \nu) (1-2 \, \nu)}}\big)\, .  \nonumber
	\end{align}
	Determining the function $\lambda_2(\lambda_1)$ using the equation $\sigma_2(\lambda_1, \lambda_2, \lambda_2) = \sigma_3(\lambda_1, \lambda_2, \lambda_2)= 0$ yields \\
	$\lambda_2(\lambda_1) = \lambda_3 (\lambda_1)= \lambda_1^{-\nu}$,
	and the tensile Cauchy stress $\sigma_1(\lambda_1)$ is given by
	\begin{equation}
		\begin{alignedat}{2}
			\label{eq:quadHenckysigma1}
			\widetilde{\sigma}_1(\lambda_1) &= \frac{1}{\lambda_1 \, \lambda_2 \, \lambda_3} \log \bigg( \big( \lambda_1^{\frac{1}{1 + \nu}} (\lambda_1 \, \lambda_2 \, \lambda_3)^{\frac{\nu}{(1+\nu) (1 - 2 \, \nu)}}\big)^{E} \bigg)\\
			&= E \cdot \frac{1}{\lambda_1 \, \lambda_1^{-\nu} \, \lambda_1^{-\nu}} \log \big( \lambda_1^{\frac{1}{1 + \nu}} (\lambda_1 \, \lambda_1^{-\nu} \, \lambda_1^{-\nu})^{\frac{\nu}{(1+\nu)(1-2 \, \nu)}}  \big) \\
			&= \frac{E}{1 + \nu} \cdot \frac{1}{\lambda_1^{1-2 \, \nu}} \log \big( \lambda_1 \, (\lambda_1^{1-2 \, \nu})^{\frac{\nu}{1- 2\, \nu}} \big)
			= \frac{E}{1 + \nu} \cdot \frac{1}{\lambda_1^{1-2 \, \nu}} \log \big( \lambda_1 \cdot\lambda_1^{\nu}\big)\\
			&= \frac{E}{1 + \nu} \cdot \frac{1}{\lambda_1^{1-2 \, \nu}} \log \big( \lambda_1^{1+\nu} \big)
			= E \cdot \lambda_1^{2 \, \nu - 1} \cdot \log \lambda_1,
		\end{alignedat}
	\end{equation}
	which, in general, is \textbf{not} monotone for $\nu \in (-1 , \frac{1}{2})$. Indeed, the incremental Young's modulus is
	\begin{equation}
		\begin{alignedat}{2}
			E^{\text{incr}} (\lambda_1) &\colonequals \frac{\dif}{\dif \lambda_1} \sigma_1 (\lambda_1)
			= E \cdot \lambda_1^{2 \, \nu - 2} ((2 \nu - 1) \cdot \log \lambda_1 + 1)
		\end{alignedat}
	\end{equation}
	and may become negative for some $\lambda_1> 1$, see Figure \ref{Exple_5_3}. 
	Regarding the energy, we have
	\begin{equation}
		\begin{alignedat}{2}
			\WW (\lambda_1, \lambda_2, \lambda_3)
			&= \mu \, ((\log \lambda_1)^2 + (\log \lambda_2)^2 + (\log \lambda_3)^2) + \frac{\lambda}{2} (\log \lambda_1 + \log \lambda_2 + \log \lambda_3)^2 \, .
		\end{alignedat}
	\end{equation}
	Inserting $\lambda_3 = \lambda_2$ and $\lambda_2(\lambda_1) = \lambda_1^{-\nu}$ we obtain
	\begin{equation}
		\begin{alignedat}{2}
			\widetilde{\WW} (\lambda_1) &= \mu ((\log \lambda_1)^2 + (\log \, (\lambda_1^{-\nu}))^2 + (\log  \, (\lambda_1^{-\nu}))^2 + \frac{\lambda}{2} (\log \lambda_1 + \log \, (\lambda_1^{-\nu}) + \log \, (\lambda_1^{-\nu}))\\
			&= \big(\mu \, (1+ 2 \, \nu^2) + \frac{\lambda}{2} (1 - 2 \, \nu)^2\big) (\log \lambda_1)^2 \, .
		\end{alignedat}
	\end{equation}
		\begin{figure}[h!]
        \hspace*{\fill}%
            \begin{minipage}[t]{0.45\textwidth}
                \centering
                \vspace{0pt}
                \includegraphics[width=\textwidth]{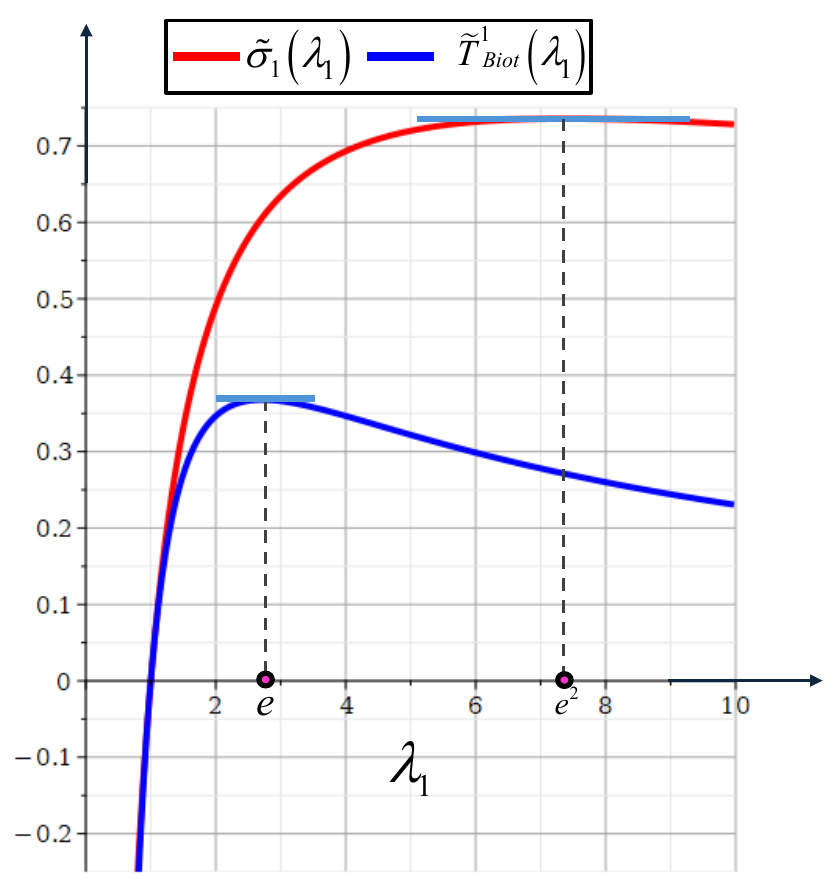}
                \captionof{figure}{Quadratic Hencky energy: tensile Cauchy stress and $T_{\Biot}^1$-stress are \textbf{non} monotone.}
                \label{Exple_5_3}
            \end{minipage}%
            \hfill
            \begin{minipage}[t]{0.45\textwidth}
                \centering
                \vspace{0pt}
                \includegraphics[width=\textwidth]{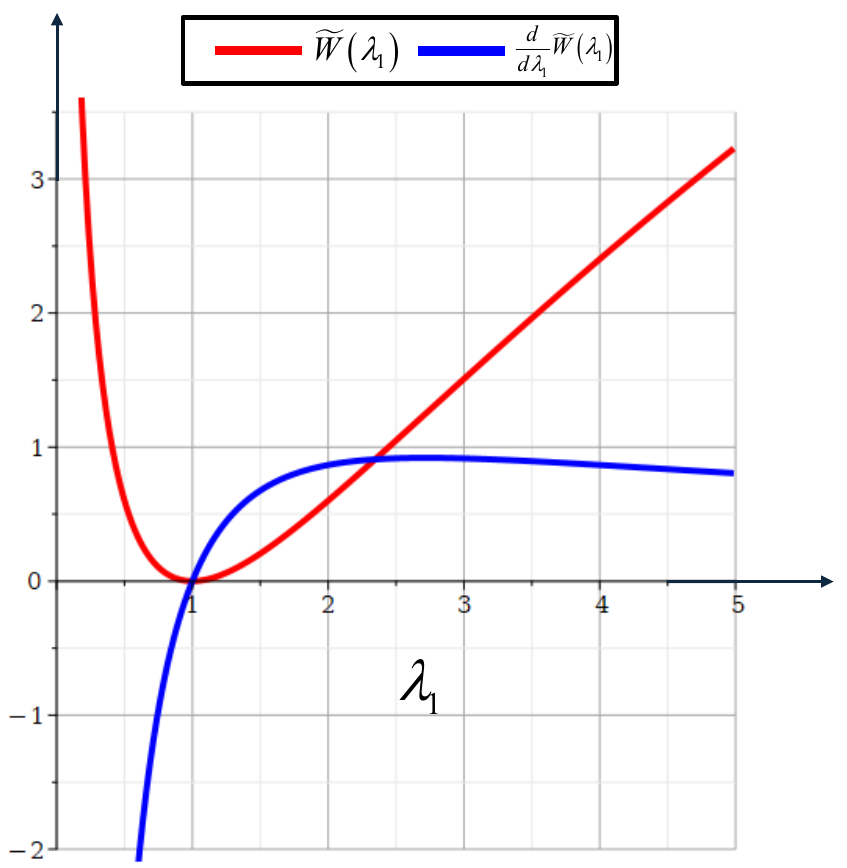}
                \captionof{figure}{Nonconvex quadratic Hencky energy in uniaxial tension and its non-monotone derivative.}
                \label{Energy_5_3}
            \end{minipage}%
        \hspace*{\fill}
        \end{figure} 
The uniaxial energy $\lambda_1 \mapsto \widetilde{\WW} (\lambda_1)$ is (slightly) not convex. Note that $\lambda_1 = \mathrm{e}, \lambda_2 = \lambda_3 = \mathrm{e}^{-\nu}$ is outside the LH-ellipticity domain (cf.~\cite{Bruhns01,agn_neff2015exponentiatedI}). Our result in \eqref{eq:compuniaxten} is consistent with the non-monotone response of $\lambda_1 \mapsto \widetilde{{\sigma}}_1 (\lambda_1)$.
\end{example}
\subsection{Uniaxial tension: incompressible response}
\label{subsec:uniaxialincompressible}
	\begin{example}[Exp-Hencky, incompressible] $\lambda_2(t) = \lambda_3(t) = \frac{1}{\sqrt{\lambda_1(t)}} = \lambda_1^{-\frac{1}{2}} (t)$, \; $\tau_2 = \tau_3 = 0$,
		\begin{equation}
			\begin{alignedat}{2}
				W_{\text{exp-Hencky}} (F) &= \frac{\mu}{k} \, \mathrm{e}^{k \norm{\log V}^2} = \frac{\mu}{k} \mathrm{e}^{k \big( (\log \lambda_1)^2 + (\log \lambda_2)^2 + (\log \lambda_3)^2\big)} \, .
			\end{alignedat}
		\end{equation}
		Choose $k=1$:
		\begin{equation}
			\widetilde{\WW}_{\text{exp-Hencky}} (\lambda_1) = \mu \, \mathrm{e}^{3/2 (\log \lambda_1)^2}
		\end{equation}
		which shows that $\lambda_1 \mapsto \widetilde{\WW}(\lambda_1)$ is convex in $\lambda_1$. Moreover
		\begin{equation}
			\begin{alignedat}{2}
				\sigma_i = \tau_i &= -p \cdot 1 + 2 \, \mu \frac{\log \lambda_i}{\lambda_i} \mathrm{e}^{k \big( (\log \lambda_1)^2 + (\log \lambda_2)^2 + (\log \lambda_3)^2 \big)}\\
			\end{alignedat}
		\end{equation}
Assuming \; $\tau_2 = \tau_3=0 $ \; determines the indeterminate pressure $p$ via 
    \begin{equation}
    	\begin{alignedat}{2}
    		p &= 2{\mkern 1mu} \mu \frac{{\log {\lambda _2}}}{{{\lambda _2}}}{{\rm{e}}^{k\left( {{{(\log {\lambda _1})}^2} + {{(\log {\lambda _2})}^2} + {{(\log {\lambda _3})}^2}} \right)}} \\
    		&= 2{\mkern 1mu} \mu \frac{{\log \lambda _1^{ - \frac{1}{2}}}}{{\lambda _1^{ - \frac{1}{2}}}}{{\rm{e}}^{k\left( {{{(\log {\lambda _1})}^2} + {{(\log (\lambda _1^{ - \frac{1}{2}}))}^2} + (\log {{(\lambda _1^{ - \frac{1}{2}})}^2}} \right)}} \\
    		&=  - {\mkern 1mu} \mu \frac{{\log {\lambda _1}}}{{\lambda _1^{ - \frac{1}{2}}}}{{\rm{e}}^{k\left( {{{(\log {\lambda _1})}^2} + 2 \cdot {{( - \frac{1}{2}\log {\lambda _1})}^2}} \right)}} =  - {\mkern 1mu} \mu \frac{{\log {\lambda _1}}}{{\lambda _1^{ - \frac{1}{2}}}}{{\rm{e}}^{k\left( {{{(\log {\lambda _1})}^2} + 2 \cdot \frac{1}{4}{{(\log {\lambda _1})}^2}} \right)}} \\
    		&=  - {\mkern 1mu} \mu \frac{{\log {\lambda _1}}}{{\lambda _1^{ - \frac{1}{2}}}}{{\rm{e}}^{\frac{3}{2}k{{(\log {\lambda _1})}^2}}} \, ,
    	\end{alignedat}
    \end{equation}
		\begin{equation}
			\begin{alignedat}{2}
				\implies \quad \widetilde{\tau}_1 &=  \mu \frac{\log \lambda_1}{\lambda_1^{-\frac{1}{2}}} \mathrm{e}^{\frac{3}{2} (\log \lambda_1)^2} + 2 \, \mu \frac{\log \lambda_1}{\lambda_1} \mathrm{e}^{\frac{3}{2} (\log \lambda_1)^2}
				= \mu \, \log \lambda_1 \, \mathrm{e}^{\frac{3}{2} (\log \lambda_1)^2} \big( \sqrt{\lambda_1} + \frac{2}{\lambda_1} \big)\\
				&= \frac{E}{2(1 + \nu)} \, \log \lambda_1 \, \mathrm{e}^{\frac{3}{2} (\log \lambda_1)^2} \big( \sqrt{\lambda_1} + \frac{2}{\lambda_1} \big)  \overset{\nu = \frac{1}{2}}{=} \frac{E}{3} \, \log \lambda_1 \, \mathrm{e}^{\frac{3}{2} (\log \lambda_1)^2} \big( \sqrt{\lambda_1} + \frac{2}{\lambda_1} \big) \\
				&\overset{\mathclap{(\lambda_1 = 1 + \delta)}}{=}\quad E \, \delta + \, \text{h.o.t.} \, .
			\end{alignedat}
		\end{equation}
The Biot stress $\widetilde{T}_{\Biot}^1 \left(\lambda_1 \right)$ follows from
\begin{eqnarray} \label{Eq5_19}
  \widetilde T_{\text{Biot}}^1\left( {{\lambda _1}} \right) = {\lambda _2} \, {\lambda _3} \, {\widetilde{\sigma} _1}\left( {{\lambda _1}} \right) = \frac{1}{{{\lambda _1}}} \widetilde{\sigma} \left( {{\lambda _1}} \right) = \frac{1}{{{\lambda _1}}} \widetilde{\tau} \left( {{\lambda _1}} \right) 
\end{eqnarray}
for incompressibility $\lambda_1 \, \lambda_2 \, \lambda_3 = 1$. The result in \eqref{eq:incunitens} implies that $\lambda_1 \mapsto \widetilde{\tau}_1 (\lambda_1)$ is monotone.
		\begin{figure}[h!]
        \hspace*{\fill}%
            \begin{minipage}[t]{0.45\textwidth}
                \centering
                \vspace{0pt}
                \includegraphics[width=\textwidth]{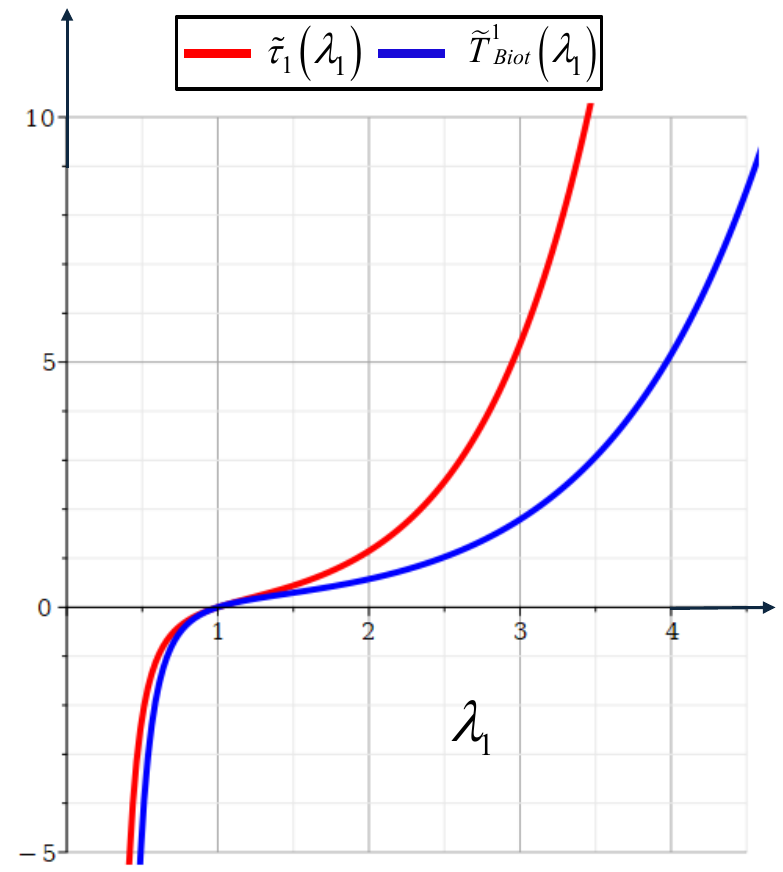}
                \captionof{figure}{Exp-Hencky incompressible: monotone tensile Kirchhoff-stress $\widetilde{\tau} _1$ and $\widetilde{T}_{\Biot}^1$-stress.}
                \label{Exple_5_4}
            \end{minipage}%
            \hfill
            \begin{minipage}[t]{0.45\textwidth}
                \centering
                \vspace{0pt}
                \includegraphics[width=\textwidth]{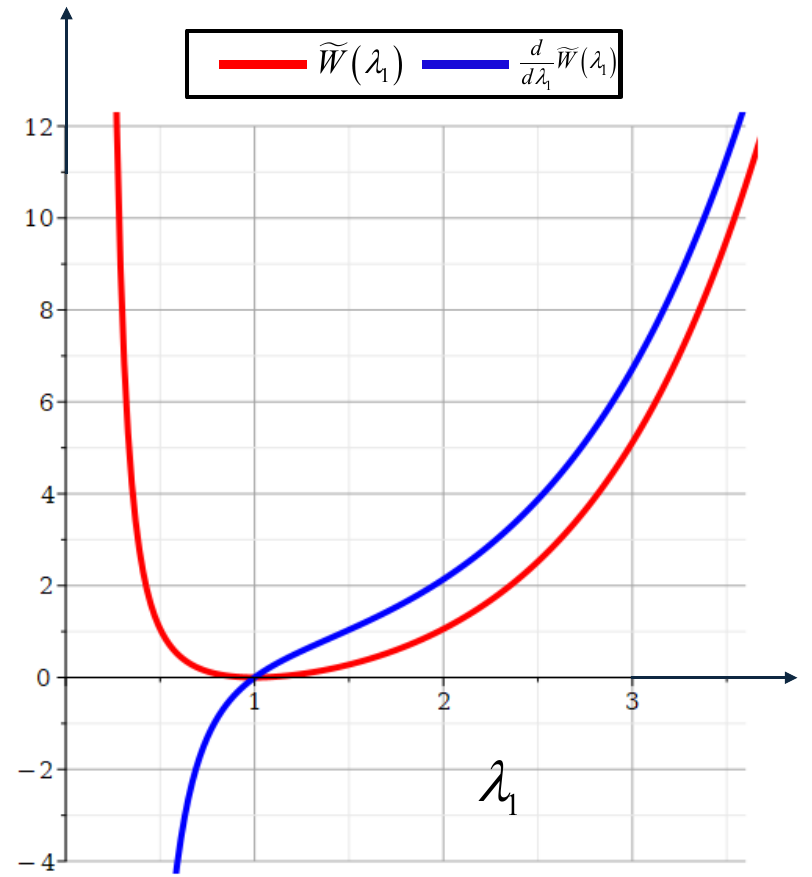}
                \captionof{figure}{Exp-Hencky incompressible: uniaxial convex energy and its monotone derivative.}
                \label{Energy_5_4}
            \end{minipage}%
        \hspace*{\fill}
        \end{figure} 
	\end{example}
\begin{example}[Neo-Hooke (cf.~\cite{Rooij2016}), incompressible]
	The Neo-Hooke model satisfies Hill's inequality in the incompressible case. We have
	\begin{eqnarray}
			\WW_{\NH}^{\nu = \frac{1}{2}} (F) = \frac{\mu}{2} \big( \norm{F}^2 - 3 \big) = \frac{\mu}{2} \left[ \tr (B) - 3 \right], \qquad \sigma_{\NH} (B) &=& \tau_{\NH} (B) = - p \, \id + \mu \, B \, .
	\end{eqnarray}
	Hill's inequality is satisfied since
	\begin{equation}
		\begin{alignedat}{2}
			\langle \tau_{\NH} (B_1) - \tau_{\NH} (B_2) , \log V_1 - \log V_2 \rangle &= \langle - p \, \id + \mu \, B_1 - (- p \, \id + \mu \, B_2) , \frac{1}{2} (\log B_1 - \log B_2) \rangle\\
			&= \frac{\mu}{2} \langle B_1 - B_2 , \log B_1 - \log B_2 \rangle > 0 \, ,
		\end{alignedat}
	\end{equation}
	due to the monotonicity of the matrix logarithm.\\
	Here again the indeterminate pressure $p$ has to be determined from the boundary conditions. In the principal stress-stretch framework this gives
	\begin{equation}
		\begin{alignedat}{2}
			\tau_i &= \sigma_i = - p \cdot 1 + \mu \, \lambda_i^2\\
			\left. \begin{matrix}
				\sigma_2 = \sigma_3 = 0 \\ 
				\tau_2 = \tau_3 = 0
			\end{matrix} \right\}
			 &\implies \quad - p + \mu \, \lambda_2^2 = 0 \iff \quad p = \mu \, \lambda_2^2 \, .
		\end{alignedat}
	\end{equation}
	Due to $\det F = 1$: $\lambda_2 = \lambda_3 = \frac{1}{\sqrt{\lambda_1}}$ and this implies
	\begin{equation}\label{Eq_23}
		\begin{alignedat}{2}
			p &= \mu \, \lambda_2^2 = \mu \big( \frac{1}{\sqrt{\lambda_1}} \big)^2 = \mu \, \frac{1}{\lambda_1}\\
			\implies \quad \widetilde{\tau}_1 (\lambda_1)&= \sigma_1 (\lambda_1) = \mu \, (\lambda_1^2 - \frac{1}{\lambda_1})
			= \frac{E}{2 \, (1+ \nu)} \, (\lambda_1^2 - \frac{1}{\lambda_1}) \overset{\nu = \frac{1}{2}}{=} \frac{E}{3} \, (\lambda_1^2 - \frac{1}{\lambda_1})\\
			&\overset{\mathclap{(\lambda_1 = 1 + \delta)}}{=} \quad E \, \delta + \, \text{h.o.t.}
		\end{alignedat}
	\end{equation}
The Biot stress $\widetilde{T}_{\text{Biot}}^1 \left(\lambda_1 \right)$ follows from above (\ref{Eq_23}) as $\widetilde T_{\text{Biot}}^1\left( {{\lambda _1}} \right) = \dfrac{1}{{{\lambda _1}}}\tilde \tau \left( {{\lambda _1}} \right)$. The result in \eqref{eq:incunitens} implies that $\lambda_1 \mapsto \widetilde{\tau}_1 (\lambda_1)$ is monotone. 
\begin{figure}[h!]
        \hspace*{\fill}%
            \begin{minipage}[t]{0.45\textwidth}
                \centering
                \vspace{0pt}
                \includegraphics[width=\textwidth]{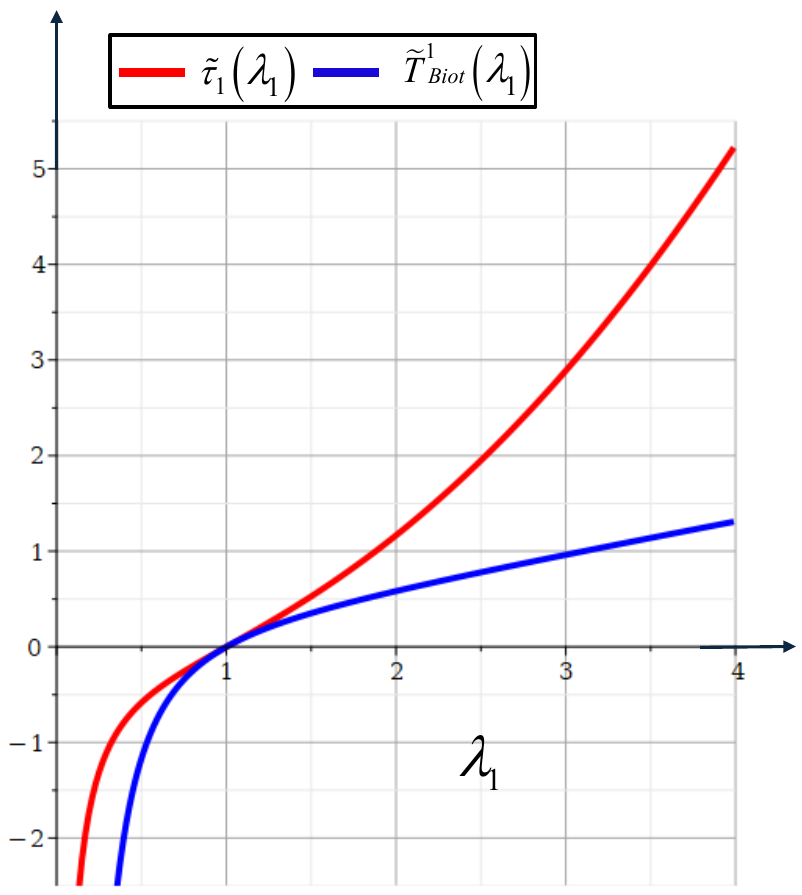}
                \captionof{figure}{Neo-Hooke incompressible: monotone tensile Kirchhoff-stress $\widetilde{\tau} _1$  and $\widetilde{T}_{\Biot}$-stress.}
                \label{Exple_5_5}
            \end{minipage}%
            \hfill
            \begin{minipage}[t]{0.45\textwidth}
                \centering
                \vspace{0pt}
                \includegraphics[width=\textwidth]{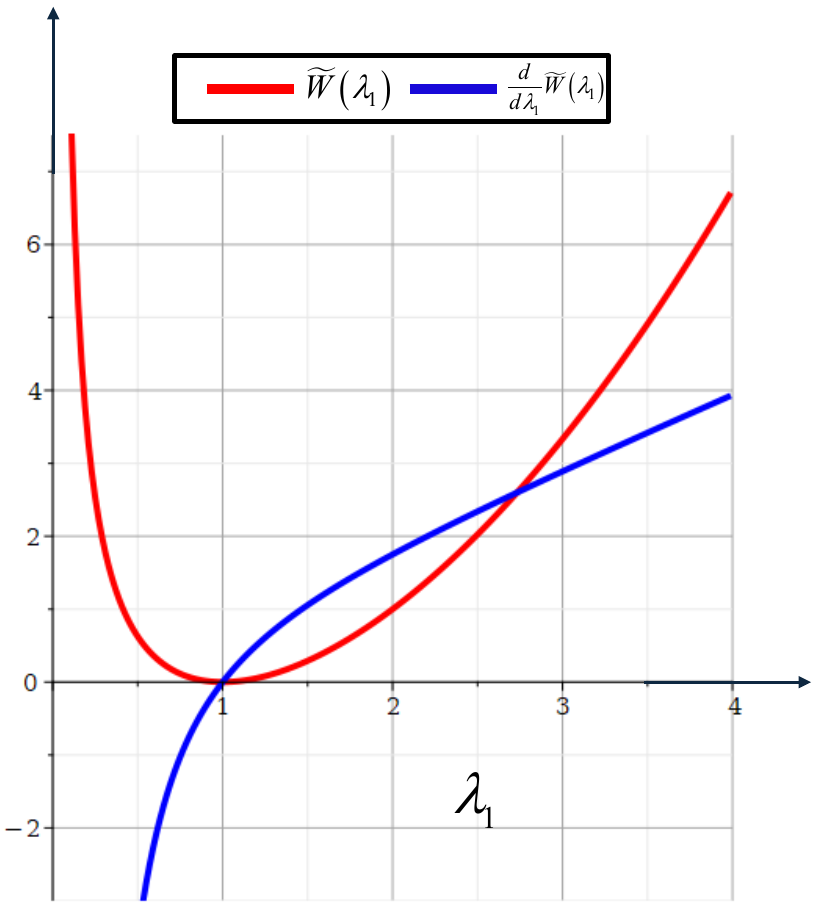}
                \captionof{figure}{Neo-Hooke incompressible: convex uniaxial energy and its monotone derivative.}
                \label{Energy_5_5}
            \end{minipage}%
        \hspace*{\fill}
        \end{figure} 
\end{example}
\begin{example}[Quadratic Hencky, incompressible]
We have then $\lambda_1(t), \lambda_2(t) = \lambda_3(t) = \frac{1}{\sqrt{\lambda_1(t)}} = \lambda_1^{-\frac{1}{2}} (t)$.\\
\begin{equation}
	\begin{alignedat}{2}
		\WW_{\text{Hencky}} &= \mu \, \norm{\log V}^2 + \frac{\lambda}{2} \, \tr^2(\log V), \\
		\widetilde{\WW}_{\text{Hencky,inc}}^{\nu = 1/2} (\lambda_1) &= \mu \big( (\log \lambda_1)^2 + (\log \lambda_2)^2 + (\log \lambda_3)^2\big)\\
		&= \mu \big( (\log \lambda_1)^2 + 2 \, (\log \lambda_2)^2 \big)
		= \mu \big( (\log \lambda_1)^2 + 2 \, (- \frac{1}{2} \log \lambda_1)^2 \big)\\
		&= \mu (1 + 2 \cdot \frac{1}{4}) (\log \lambda_1)^2
		= \frac{3}{2} \, \mu \, (\log \lambda_1)^2
		= \frac{3}{2} \frac{E}{2(1 + \nu)}\big\vert_{\nu = \frac{1}{2}} (\log \lambda_1)^2 
		= \frac{E}{2} (\log \lambda_1)^2 \, .
	\end{alignedat}
\end{equation}
Note that $\lambda_1 \mapsto \widetilde{\WW}_{\text{Hencky,inc}}^{\nu = 1/2} (\lambda_1)$ is \textbf{not} convex in $\lambda_1$. 
	
With \eqref{eq:quadHenckysigma1} we obtain for $\nu = \frac{1}{2}$: $\widetilde{\sigma} _1 (\lambda_1) = \widetilde{\tau} _1 (\lambda_1) = E \, \log \lambda_1$.
A direct computation yields as well
\begin{equation}
	\begin{alignedat}{3}
		\tau_i &= \sigma_i = - p \cdot 1 + 2 \, \mu \log \lambda_i \, , \quad
		\tau_2 = \tau_3 = 0,\\
		 \implies p &= 2 \, \mu \, \log \lambda_2 = 2 \, \mu (- \frac{1}{2} \log \lambda_1)
		= - \mu \, \log \lambda_1\\
		\implies \widetilde{\tau}_1 &= 2 \, \mu \log \lambda_1 + \mu \log \lambda_1 =3 \,  \mu \, \log \lambda_1 = 3 \, \frac{E}{2 \, (1 + \nu)} \log \lambda_1 \overset{\nu = \frac{1}{2}}{=} E \, \log \lambda_1
		\overset{(\lambda_1 = 1 + \delta)}{=} E \, \delta + \, \text{h.o.t.}
	\end{alignedat}
\end{equation}
We observe that $\lambda_1 \mapsto \widetilde{\tau} _1(\lambda_1)$ is monotone as it should be due to the satisfaction of Hill's inequality in the incompressible case and our statement in \eqref{eq:incunitens}.
		\begin{figure}[h!]
        \hspace*{\fill}%
            \begin{minipage}[t]{0.45\textwidth}
                \centering
                \vspace{0pt}
                \includegraphics[width=\textwidth]{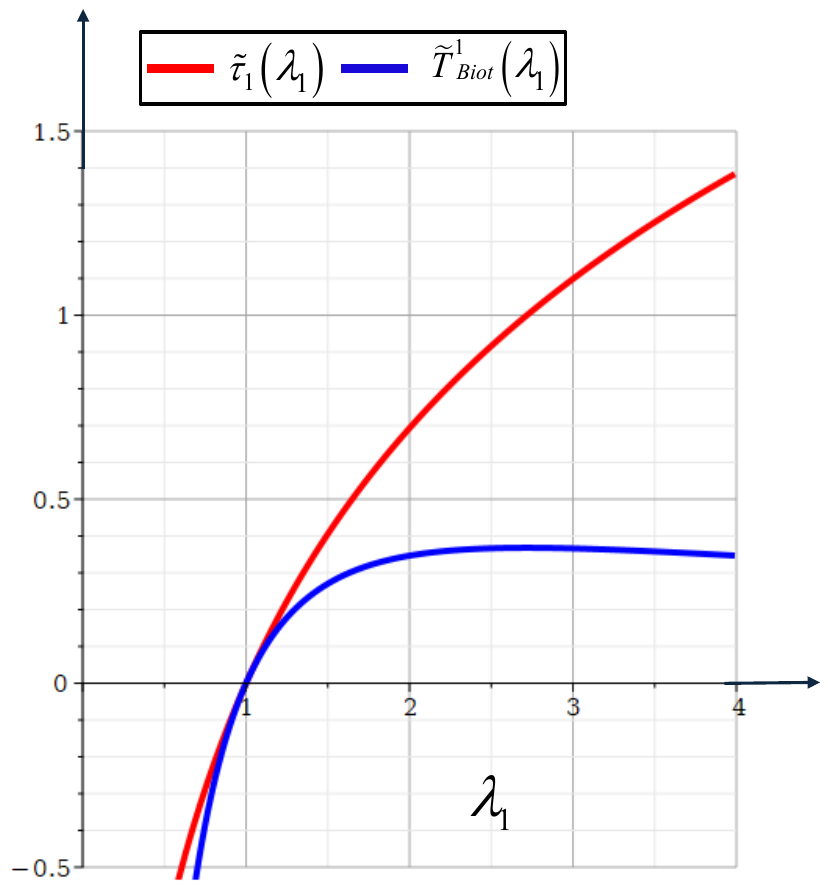}
                \captionof{figure}{Quadratic Hencky incompressible: the uniaxial Kirchhoff-stress $\tau_1$ is monotone, while the $T_{\Biot}^1$-stress remains non monotone.}
                \label{Exple_5_6}
            \end{minipage}%
            \hfill
            \begin{minipage}[t]{0.45\textwidth}
                \centering
                \vspace{0pt}
                \includegraphics[width=\textwidth]{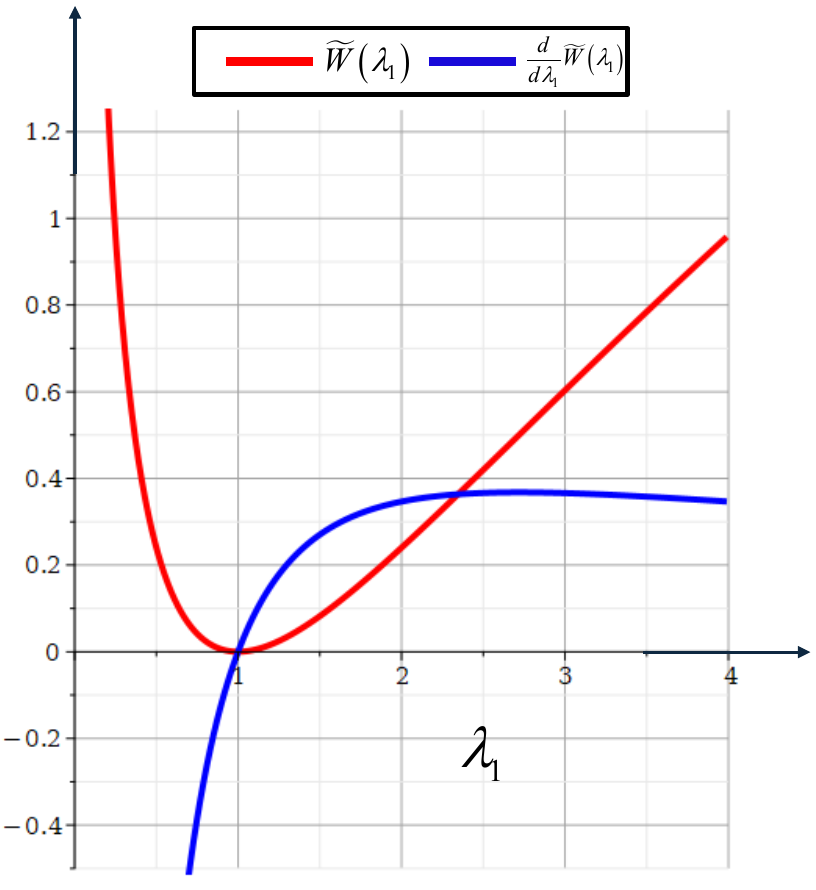}
                \captionof{figure}{Quadratic Hencky incompressible: the uniaxial energy remains non convex and its derivative is non monotone.}
                \label{Energy_5_6}
            \end{minipage}%
        \hspace*{\fill}
        \end{figure} 
\end{example}
\section{Conclusion}
We have clarified that (CSP) as constitutive assumption is different from a notion of positive second order internal work in nonlinear elasticity while it coincides formally with the Drucker stability postulate for geometrically linear kinematics. Following, we have shown that (CSP) $\iff$ (TSTS-M$^{++}$) simplifies considerably for a special family of equilibrium solutions: namely if $t \mapsto F(t)$ is homogeneous and diagonal. In this case, the corotational derivative $\frac{\DD^{\circ}}{\DD t}$ reduces to the material time derivative $\frac{\DD}{\DD t}$ ($F$ is diagonal) and the material time derivative $\frac{\DD}{\DD t}$ reduces to the usual partial time derivative $\partial_t$ for homogeneous solutions. Applying this insight to standard experimental tests like uniaxial tension, equibiaxial extension, planar tension and hydrostatic tension shows that the corresponding Cauchy stress incremental moduli are always positive if (CSP) is assumed. Three examples for uniaxial tension with and without (CSP) are explicitly worked out, for the compressible and the incompressible case, Cauchy (Kirchhoff). For incompressibility, Hill's inequality coincides with (CSP) and already implies positive incremental moduli.\\

Thus, (CSP) $\iff$ (TSTS-M$^{++}$) emerges as a suitable minimal \textbf{constitutive stability assumption} in isotropic nonlinear elasticity, complementing \textbf{local material stability}, here considered as satisfaction of the LH-ellipticity condition
\begin{equation}
	\DD_F^2 \WW (F) . (\xi \otimes \eta, \xi \otimes \eta) > 0 \quad \forall \ \xi, \eta \neq 0
\end{equation}
which implies stability of the homogeneous state $F = \DD \varphi$ against infinitesimal interior perturbations (cf.~\cite{vanHove1947}). For a stable idealized elastic material with physically reasonable response, it remains therefore to find an isotropic hyperelastic formulation that satisfies simultaneously (CSP) and LH-ellipticity throughout.\\

\footnotesize
\section*{References}
\printbibliography[heading=none]

@article{agn_neff2014logarithmic,
	title={A logarithmic minimization property of the unitary polar factor in the spectral and Frobenius norms},
	author={Neff, Patrizio and Nakatsukasa, Yuji and Fischle, Andreas},
	journal={SIAM Journal on Matrix Analysis and Applications},
	volume={35},
	number={3},
	pages={1132--1154},
	year={2014},
	publisher={SIAM},
	doi={10.1137/130909949},
	class =	{matrix}
}

@article{agn_lankeit2014minimization,
	title={The minimization of matrix logarithms: On a fundamental property of the unitary polar factor},
	author={Lankeit, Johannes and Neff, Patrizio and Nakatsukasa, Yuji},
	journal={Linear Algebra and its Applications},
	volume={449},
	pages={28--42},
	year={2014},
	publisher={Elsevier},
	doi={10.1016/j.laa.2014.02.012},
	class =	{matrix}
}

@article{agn_neff2014grioli,
	title={On Grioli's minimum property and its relation to {C}auchy's polar decomposition},
	author={Neff, Patrizio and Lankeit, Johannes and Madeo, Angela},
	journal={International Journal of Engineering Science},
	volume={80},
	pages={209--217},
	year={2014},
	publisher={Elsevier},
	doi={10.1016/j.ijengsci.2014.02.026},
	class =	{matrix}
}

@article{agn_neff2015exponentiatedI,
	title={The exponentiated {H}encky-logarithmic strain energy. {P}art {I}: {C}onstitutive issues and rank-one convexity},
	author={Neff, Patrizio and Ghiba, Ionel-Dumitrel and Lankeit, Johannes},
	journal={Journal of Elasticity},
	volume={121},
	number={2},
	pages={143--234},
	year={2015},
	publisher={Springer},
	doi={10.1007/s10659-015-9524-7},
	class =	{mathematics}
}

@article{agn_neff2015exponentiatedII,
	title={The exponentiated Hencky-logarithmic strain energy. {P}art {II}: coercivity, planar polyconvexity and existence of minimizers},
	author={Neff, Patrizio and Lankeit, Johannes and Ghiba, Ionel-Dumitrel and Martin, Robert J. and Steigmann, David J.},
	journal={Zeit\-schrift f{\"u}r angewandte Mathematik und Physik},
	volume={66},
	number={4},
	pages={1671--1693},
	year={2015},
	publisher={Springer},
	doi={10.1007/s00033-015-0495-0},
	class =	{mathematics}
}

@article{agn_ghiba2015ellipticity,
	title={An ellipticity domain for the distortional Hencky logarithmic strain energy},
	author={Ghiba, Ionel-Dumitrel and Neff, Patrizio and Martin, Robert J.},
	journal={Proceedings of the Royal Society of London A: Mathematical and Physical Sciences},
	volume={471},
	number={2184},
	year={2015},
	publisher={The Royal Society},
	doi={10.1098/rspa.2015.0510},
	class =	{mathematics}
}

@article{agn_neff2015geometry,
	title =		{Geometry of logarithmic strain measures in solid mechanics},
	author =	{Neff, Patrizio and Eidel, Bernhard and Martin, Robert J.},
	journal =	{Archive for Rational Mechanics and Analysis},
	volume =	{222},
	number =	{2},
	pages =		{507--572},
	note =		{\availableatarxiv{1505.02203}},
	year =		{2016},
	doi =		{10.1007/s00205-016-1007-x},
	class =	{mathematics}
}

@article{agn_neff2014axiomatic,
	title={The axiomatic deduction of the quadratic Hencky strain energy by Heinrich Hencky},
	author={Neff, Patrizio and Eidel, Bernhard and Martin, Robert J},
	journal={\arxivjournal{1402.4027}},
	year={2014},
	class =	{mathematics}
}

@article{agn_martin2018non,
	title={A non-ellipticity result, or the impossible taming of the logarithmic strain measure},
	author={Martin, Robert J and Ghiba, Ionel-Dumitrel and Neff, Patrizio},
	journal={International Journal of Non-Linear Mechanics},
	volume={102},
	pages={147--158},
	year={2018},
	publisher={Elsevier},
	class =	{mathematics}
}

@article{agn_nedjar2018finite,
	title={A finite element implementation of the isotropic exponentiated Hencky-logarithmic model and simulation of the eversion of elastic tubes},
	author={Nedjar, Boumediene and Baaser, Herbert and Martin, Robert J and Neff, Patrizio},
	journal={Computational Mechanics},
	volume={62},
	number={4},
	pages={635--654},
	year={2018},
	publisher={Springer},
	 note =		{\availableatarxiv{1705.08381}},
	 class =	{numerics}
}

@article{jog2013,
  year =			{2013},
  issn =			{0939-1533},
  journal =		{Archive of Applied Mechanics},
  volume =		{83},
  number =		{5},
  doi =			{10.1007/s00419-012-0711-8},
  title =		{Conditions for the onset of elastic and material instabilities in hyperelastic materials},
  publisher =	{Springer-Verlag},
  author =		{Jog, C. S. and Patil, K. D.},
  pages =		{661-684}
}

@article{richter1948,
  year =			{1948},
  journal =		{Zeitschrift f{\"u}r Angewandte Mathematik und Mechanik},
  volume =		{28},
  number =		{7/8},
  title =		{Das isotrope {E}lastizit{\"a}tsgesetz},
  author =		{Richter, H.},
  pages =		{205--209},
  note =		{\availableaturl[uni-due.de]{https://www.uni-due.de/imperia/md/content/mathematik/ag_neff/richter_isotrop_log.pdf}}
}

@article{bakerEri54,
  author =	{M. Baker and J. L. Ericksen},
  title =	{Inequalities restricting the form of the stress-deformation relation for isotropic elastic solids and {Reiner}-{Rivlin} fluids},
  journal =	{Journal of the Washington Academy of Sciences},
  year =	{1954},
  volume =	{44},
  pages =	{33--35}
}

@article{Bruhns01,
  author =	{O. T. Bruhns and H. Xiao and A. Mayers},
  title =	{Constitutive inequalities for an isotropic elastic strain energy function based on {H}encky's logarithmic strain tensor},
  journal =	{Proceedings of the Royal Society of London A: Mathematical and Physical Sciences},
  volume=	{457},
  pages=	{2207-2226},
  year =	{2001}
}

@book{marsden1994foundations,
  title =		{Mathematical Foundations of Elasticity},
  author =		{Marsden, J. E. and Hughes, T.},
  year =		{1994},
  publisher =	{Courier Dover Publications}
}

@article{doyle1956nonlinear,
  title =		{Nonlinear elasticity},
  author =		{Doyle, T. C. and Ericksen, J. L.},
  journal =		{Advances in Applied Mechanics},
  volume =		{4},
  pages =		{53--115},
  year =		{1956},
  publisher =	{Elsevier}
}

@article{CSP2024,
author	= {P. Neff and S. Holthausen and M. V. d'Agostino and D. Bernardini and A. Sky and I. D. Ghiba and R. J. Martin},
title		= {Hypo-elasticity, {C}auchy-elasticity, corotational stability and monotonicity in the logarithmic strain},
journal	= {submitted, arXiv:2409.20051},
year={2024}
}

@article{sidoroff1974restrictions,
	author	= {R. Sidoroff.},
	title		= {Sur les restrictions \`{a} imposer \`{a} l'\'{e}nergie de d\'{e}formation d'un mat\'{e}riau hyper\'{e}lastique.},
	journal	= {Comptes Rendus de l'Académie des Sciences Paris},
	volume	= {279},
	pages		= {379--382},
	year		= {1974}
}

@article{mandel1966,
author = {J. Mandel.},
title = {Conditions de stabilité et postulat de Drucker},
journal = {In: J. Kravtchenko, P. M. Sirieys (eds) Rheology and Soil Mechanics/Rhéologie et Mécanique des Sols, International Union of Theoretical and Applied Mechanics},
pages = {58--68},
year = {1966},
publisher = {Springer}
}

@article{hill1957,
author = {R. Hill},
title = {On uniqueness and stability in the theory of finite elastic strain},
journal = {Journal of the Mechanics and Physics of Solids},
volume = {5},
number = {4},
pages = {229--241},
year = {1957}
}

@article{hill1968constitutive,
author = {R. Hill},
title = {On constitutive inequalities for simple materials\,-\, {I}.},
journal = {Journal of the Mechanics and Physics of Solids},
volume = {16},
number = {4},
year = {1968}
}

@article{hill1970constitutive,
author = {R. Hill},
title = {Constitutive inequalities for isotropic elastic solids under finite strain},
journal = {Proceedings of the Royal Society A: Mathematical, Physical and Engineering Sciences},
volume = {314},	
number = {1519},
pages = {457--472},
year = {1970}
}

@article{leblond1992constitutive,
author = {J. B. Leblond},
title = {A constitutive inequality for hyperelastic materials in finite strain},
journal = {European Journal of Mechanics - A/Solids},
volume = {11},
number = {4},
pages = {447--466},
year = {1992}
}

@article{Leblond2024,
author = {M. V. d'Agostino and S. Holthausen and D. Bernardini and A. Sky and I. D. Ghiba and R. J. Martin and P. Neff},
title = {A constitutive condition for idealized isotropic Cauchy elasticity involving the logarithmic strain},
journal = {to appear in Journal of Elasticity, arXiv: 2409.01811},
year = {2024}
}

@article{poscor2024,
author = {P. Neff and S. Holthausen and S. N. Korobeynikov and I. D. Ghiba and R. J. Martin},
title = {A natural requirement for objective corotational rates - on structure preserving corotational rates},
journal = {to appear in Acta Mechanica, arXiv: 2409.19707},
year = {2024}
}

@article{vanHove1947,
author = {L. C. P. van Hove},
title = {Sur l'extension de la condition de Legendre du calcul des variations aux intégrales multiples a plusieurs fonctions inconnues},
journal = {Proceedings of the Koninklijke Nederlandse Akademie van Wetenschappen},
volume = {50},
number = {1},
pages = {18--23},
year = {1947}
}

@article{richter1949hauptaufsatze,
author = {H. Richter},
title = {Verzerrungstensor, Verzerrungsdeviator und Spannungstensor bei endlichen Formänderungen},
journal = {Zeitschrift für Angewandte Mathematik und Mechanik},
volume = {29},
number = {3},
pages = {65--75},
year = {1949}
}

@article{petryk1985,
author = {H. Petryk},
title = {On the second-order work in plasticity},
journal = {Archives of Mechanics (Warszawa)},
volume = {43},
number = {2-3},
pages = {377--397},
year = {1991}
}

@article{mihai2017,
author = {L. A. Mihai and A. Goriely},
title = {How to characterize a nonlinear elastic material? A review on nonlinear constitutive parameters in isotropic finite elasticity},
journal ={Proceedings of the Royal Society A: Mathematical, Physical and Engineering Sciences},
volume = {473},
number = {2207},
pages = {20170607},
year = {2017}
}

@article{majorsym2024,
author = {S. Federico and S. Holthausen and N. J. Husemann and P. Neff},
title = {Major symmetry of the induced tangent stiffness tensor for the Zaremba-Jaumann rate and Kirchhoff stress in hyperelasticity: two different approaches},
journal = {submitted, arXiv: 2410.22163},
year = {2024}
}

@article{yavari2024,
author = {A. Yavari},
title = {Universal deformations and inhomogeneities in isotropic Cauchy elasticity},
journal = {arXiv:2404.06235},
year = {2024}
}

@article{scott2006,
author = {N. H. Scott},
title = {The incremental bulk modulus, Young's modulus and Poisson's ratio in nonlinear isotropic elasticity: physically reasonable response},
journal = {Mathematics and Mechanics of Solids},
volume = {12},
pages = {526--542},
year = {2006}
}

@article{Drucker1950,
author = {D. C. Drucker},
title = {Some implications of work hardening and ideal plasticity},
journal = {Quarterly of Applied Mathematics},
volume = {7},
number = {4},
pages = {411--418},
year = {1950}
}

@article{Ericksen1955,
author = {J. L. Ericksen},
title = {Deformations possible in every compressible, isotropic, perfectly elastic material},
journal = {Studies in Applied Mathematics},
volume = {34},
number = {1-4},
pages = {126--128},
year = {1955}
}

@article{Young2024,
author = {P. Neff and N. J. Husemann and S. Holthausen and A. S. Nguetcho Tchakoutio and I. D. Ghiba and R. J. Martin},
title = {An essay on constitutive stability in idealized isotropic nonlinear elasticity for universal deformations, positive incremental moduli and the onset of necking},
journal = {in preparation}
}

@article{Martin2024,
author = {R. J. Martin and I. D. Ghiba and P. Neff},
title = {The corotational stability postulate is equivalent to the true stress-true strain monotonicity condition},
journal = {in preparation}
}

@article{Rooij2016,
author = {R. de Rooij and E. Kuhl},
title = {Constitutive modeling of brain tissue: current perspective},
journal = {Applied Mechanics Review},
volume = {68},
pages = {010801},
year = {2016}
}

@article{Ghiba2024,
author = {I. D. Ghiba and R. J. Martin and P. Neff},
title = {Constitutive properties for isotropic energies in ideal nonlinear elasticity for solid materials: numerical evidence for invertibility and monotonicity in different stress-strain pairs},
journal = {in preparation}
}

@article{Martin2024Hill,
author = {R. J. Martin and J. Voss and I. D. Ghiba and M. V. d'Agostino and P. Neff},
title = {Monotonicity of isotropic tensor functions on the set of symmetric matrices: Hill's generalization of the Davis-Lewis convexity theorem revised},
journal = {in preparation}
}

@article{Drucker1951,
author = {D. C. Drucker},
title = {A more fundamental approach to plastic stress-strain relations},
journal = {Proceedings of the First U. S. National Congress of Applied Mechanics},
pages = {187--491},
year = {1951}
}

@article{Hill1958,
author = {R. Hill},
title = {A general theory of uniqueness and stability in elastic-plastic solids},
journal = {Journal of the Mechanics and Physics of Solids},
volume ={6},
number = {3},
pages = {236--249},
year = {1958}
}

@article{Hill1959,
author = {R. Hill},
title = {Some basic principles in the mechanics of solids without a natural time},
journal = {Journal of the Mechanics and Physics of Solids},
volume = {7},
pages = {209--225},
year = {1959}
}

@article{biezeno1928,
	author	= {C. B. Biezeno and H. Hencky},
	title		= {On the general theory of elastic stability.},
	journal	= {Koninklijke Akademie van Wettenschappen te Amsterdam},
	volume	= {31},
	pages		= {569--592},
	year		= {1928},
}

@article{korobeynikov2023,
	author	= {S. N. Korobeynikov.},
	title		= {Families of {H}ooke-like isotropic hyperelastic material models and their rate formulations.},
	journal	= {Archive of Applied Mechanics},
	volume	= {93},
	pages		= {3863--3893},
	year		= {2023},
}

@article{thiel2018,
author = {C. Thiel, J. Voss, R. J. Martin and P. Neff},
title = {Shear, pure and simple},
journal = {International Journal of Non-Linear Mechanics, arXiv: 1806.07749},
volume = {112},
pages = {57--72},
year = {2018/9}
}
\normalsize
\begin{appendix}
\section{Notation}
\label{sec:app_notation}
\begin{figure}[h!]
	\begin{center}		
		\begin{minipage}[h!]{0.8\linewidth}
			\centering
			\hspace*{-40pt}
			\includegraphics[scale=0.4]{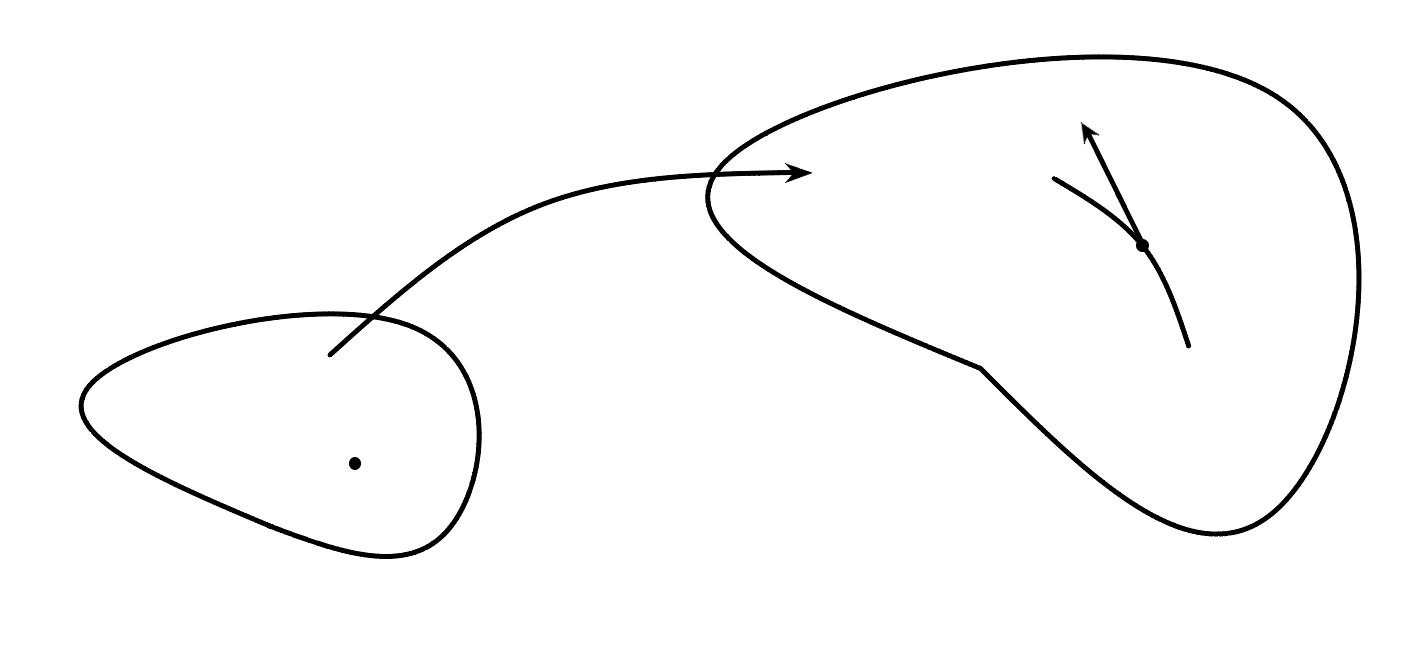}
			\put(-40,30){\footnotesize $\Omega_\xi$}
			\put(-340,25){\footnotesize $\Omega_x$}
			\put(-316,64){\footnotesize $x$}
			\put(-280,148){\footnotesize $\varphi(x,t)$}
			\put(-104,168){\footnotesize $\overline V(x,t) \!=\! v(\xi,t)$}
			\put(-88,119){\footnotesize $\xi$}
			\put(-105,90){\footnotesize curve $t \mapsto \varphi(x,t)$}
			\put(-85,80){\footnotesize  for $x$ fixed}
		\end{minipage} 
		\caption{Illustration of the deformation $\varphi(x,t): \Omega_x \to \Omega_{\xi}$ and the velocity $\overline V(x,t) = v(\xi,t)$.}
		\label{yfig1}
	\end{center}
\end{figure}
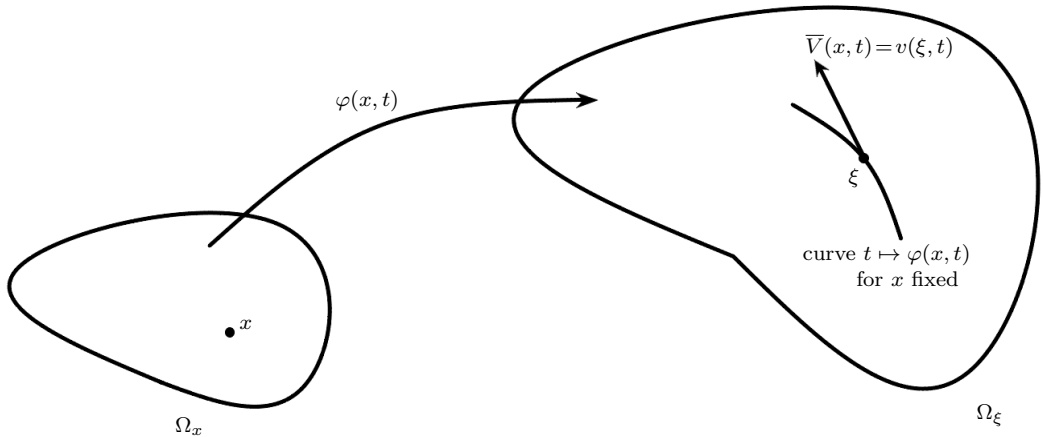

\textbf{The deformation $\varphi(x,t)$, the material time derivative $\frac{\DD}{\DD t}$ and the partial time derivative $\partial_t$} \\
\\
In accordance with \cite{marsden1994foundations} we agree on the following convention regarding an elastic deformation $\varphi$ and time derivatives of material quantities:

Given two sets $\Omega, \Omega_{\xi} \subset \R^3$ we denote by $\varphi: \Omega \to \Omega_{\xi}, x \mapsto \varphi(x) = \xi$ the deformation from the \emph{reference configuration} $\Omega$ to the \emph{current configuration} $\Omega_{\xi}$. A \emph{motion} of $\Omega$ is a time-dependent family of deformations, written $\xi = \varphi(x,t)$. The \emph{velocity} of the point $x \in \Omega$ is defined by $\overline{V}(x,t) = \partial_t \varphi(x,t)$ and describes a vector emanating from the point $\xi = \varphi(x,t)$ (see also Figure \ref{yfig1}). Similarly, the velocity viewed as a function of $\xi \in \Omega_{\xi}$ is denoted by $v(\xi,t)$. 

Considering an arbitrary material quantity $Q(x,t)$ on $\Omega$, equivalently represented by $q(\xi,t)$ on $\Omega_\xi$, we obtain by the chain rule for the time derivative of $Q(x,t)$
\begin{align}
	\frac{\DD}{\DD t}q(\xi,t) \colonequals \frac{\dif}{\dif t}[Q(x,t)] = \DD_\xi q(\xi,t).v(\xi,t) + \partial_t q(\xi,t) \, .
\end{align}
Since it is always possible to view any material quantity $Q(x,t) = q(\xi,t)$ from two different angles, namely by holding $x$ or $\xi$ fixed, we agree to write
\begin{itemize}
	\item $\dot q \colonequals \dd \frac{\DD}{\DD t}[q]$ for the material (substantial) derivative of $q$ with respect to $t$ holding $x$ fixed and
	\item $\partial_t q$ for the derivative of $q$ with respect to $t$ holding $\xi$ fixed.
\end{itemize}
For example, we obtain the velocity gradient $L := \DD_\xi v(\xi,t)$ by
\begin{align}
	L = \DD_\xi v(\xi,t) = \DD_\xi \overline V(x,t) &\overset{\text{def}}{=} \DD_\xi \frac{\dif}{\dif t} \varphi(x,t) = \DD_{\xi} \partial_t \varphi(\varphi^{-1}(\xi,t),t) = \partial_t \DD_x \varphi(\varphi^{-1}(\xi,t),t) \, \DD_\xi \big(\varphi^{-1}(\xi,t)\big) \notag \\
	&=  \partial_t \DD_x \varphi(\varphi^{-1}(\xi,t),t) \, (\DD_x \varphi)^{-1}(\varphi^{-1}(\xi,t),t) = \dot F(x,t) \, F^{-1}(x,t) = L \, ,
\end{align}
where we used that $\partial_t = \frac{\dif}{\dif t} = \frac{\DD}{\DD t}$ are all the same, if $x$ is fixed. \\
\\
As another example, when determining a corotational rate $\frac{\DD^{\circ}}{\DD t}$ we write
\begin{align}
	\frac{\DD^{\circ}}{\DD t}[\sigma] = \frac{\DD}{\DD t}[\sigma] + \sigma \, \Omega^{\circ} - \Omega^{\circ} \, \sigma = \dot \sigma + \sigma \, \Omega^{\circ} - \Omega^{\circ} \, \sigma \, .
\end{align}
However, if we solely work on the current configuration, i.e.~holding $\xi$ fixed, we write $\partial_t v$ for the time-derivative of the velocity (or any quantity in general). \\
\\
\noindent \textbf{Inner product} \\
\\
For $a,b\in\R^n$ we let $\langle {a},{b}\rangle_{\R^n}$  denote the scalar product on $\R^n$ with associated vector norm $\norm{a}_{\R^n}^2=\langle {a},{a}\rangle_{\R^n}$. We denote by $\R^{n\times n}$ the set of real $n\times n$ second order tensors, written with capital letters. The standard Euclidean scalar product on $\R^{n\times n}$ is given by
$\langle {X},{Y}\rangle_{\R^{n\times n}}=\tr{(X Y^T)}$, where the superscript $^T$ is used to denote transposition. Thus the Frobenius tensor norm is $\norm{X}^2=\langle {X},{X}\rangle_{\R^{n\times n}}$, where we usually omit the subscript $\R^{n\times n}$ in writing the Frobenius tensor norm. The identity tensor on $\R^{n\times n}$ will be denoted by $\id$, so that $\tr{(X)}=\langle {X},{\id}\rangle$. \\
\\
\noindent \textbf{Frequently used spaces} 
\begin{itemize}
	\item $\Sym(n), \rm \Sym^+(n)$ and $\Sym^{++}(n)$ denote the symmetric, positive semi-definite symmetric and positive definite symmetric second order tensors respectively.
	\item ${\rm GL}(n):=\{X\in\R^{n\times n}\;|\det{X}\neq 0\}$ denotes the general linear group.
	\item ${\rm GL}^+(n):=\{X\in\R^{n\times n}\;|\det{X}>0\}$ is the group of invertible matrices with positive determinant.
	\item $\mathrm{O}(n):=\{X\in {\rm GL}(n)\;|\;X^TX=\id\}$.
	\item ${\rm SO}(n):=\{X\in {\rm GL}(n,\R)\;|\; X^T X=\id,\;\det{X}=1\}$.
	\item $\mathfrak{so}(3):=\{X\in\mathbb{R}^{3\times3}\;|\;X^T=-X\}$ is the Lie-algebra of skew symmetric tensors.
	\item The set of positive real numbers is denoted by $\R_+:=(0,\infty)$, while $\overline{\R}_+=\R_+\cup \{\infty\}$.
\end{itemize}
\textbf{Frequently used tensors}
\begin{itemize}
	\item $F = \DD \varphi(x,t)$ is the Fréchet derivative (Jacobian matrix) of the deformation $\varphi(\cdot,t) : \Omega_x \to \Omega_{\xi} \subset \R^3$. $\varphi(x,t)$ is usually assumed to be a diffeomorphism at every time $t \ge 0$ so that the inverse mapping $\varphi^{-1}(\cdot,t) : \Omega_{\xi} \to \Omega_x$ exists.
	\item $C=F^T \, F$ is the right Cauchy-Green strain tensor.
	\item $B=F\, F^T$ is the left Cauchy-Green (or Finger) strain tensor.
	\item $U = \sqrt{F^T \, F} \in \Sym^{++}(3)$ is the right stretch tensor, i.e.~the unique element of ${\rm Sym}^{++}(3)$ with $U^2=C$.
	\item $V = \sqrt{F \, F^T} \in \Sym^{++}(3)$ is the left stretch tensor, i.e.~the unique element of ${\rm Sym}^{++}(3)$ with $V^2=B$.
	\item $\log V = \frac12 \, \log B$ is the spatial logarithmic strain tensor or Hencky strain.
	\item We write $V = Q$ diag($\lambda_1, \lambda_2, \lambda_3$) $Q^T$, where $\lambda_i \in \R_+$ are the principal stretches.
	\item $L = \dot F \, F^{-1} = \DD_\xi v(\xi)$ is the spatial velocity gradient.
	\item $D = \sym \, L$ is the spatial rate of deformation, the Eulerian strain rate tensor.
	\item $W = \sk \, L$ is the vorticity tensor.
	\item We also have the polar decomposition $F = R \, U = V R \in {\rm GL}^+(3)$ with an orthogonal matrix $R \in \OO(3)$ (cf. Neff et al.~\cite{agn_neff2014grioli}), see also \cite{agn_lankeit2014minimization,agn_neff2014logarithmic}.
\end{itemize}
\noindent \textbf{Tensor domains} \\
\\
Denoting the reference configuration by $\Omega_x$ with tangential space $T_x \Omega_x$ and the current/spatial configuration by $\Omega_\xi$ with tangential space $T_\xi \Omega_\xi$ as well as $\varphi(x) = \xi$, we have the following relations (see also Figure \ref{yfig2}):

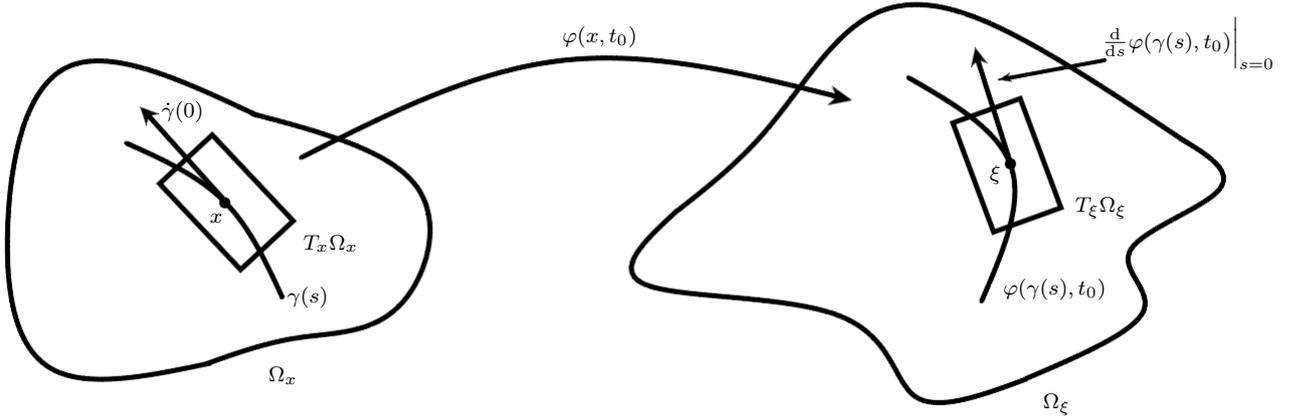
\begin{figure}[h!]
	\begin{center}		
		\begin{minipage}[h!]{0.8\linewidth}
			\centering
			\hspace*{-80pt}
			\includegraphics[scale=0.5]{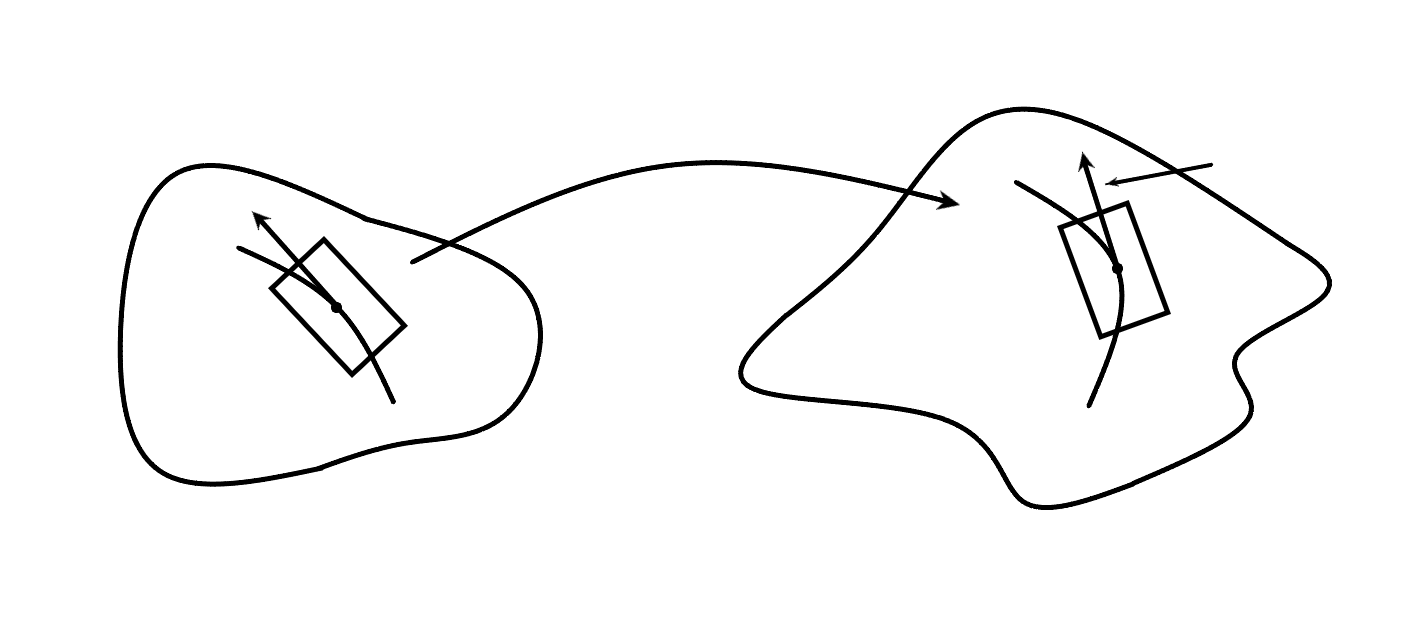}
			\put(-100,45){\footnotesize $\Omega_\xi$}
			\put(-390,55){\footnotesize $\Omega_x$}
			\put(-412,115){\footnotesize $x$}
			\put(-430,155){\footnotesize $\dot \gamma(0)$}
			\put(-377,105){\footnotesize $T_x \Omega_x$}
			\put(-383,85){\footnotesize $\gamma(s)$}
			\put(-280,183){\footnotesize $\varphi(x,t_0)$}
			\put(-120,131){\footnotesize $\xi$}
			\put(-78,181){\footnotesize $\frac{\dif}{\dif s}\varphi(\gamma(s),t_0)\bigg\vert_{s=0}$}
			\put(-88,119){\footnotesize $T_\xi \Omega_\xi$}
			\put(-115,88){\footnotesize $\varphi(\gamma(s),t_0)$}
		\end{minipage} 
		\caption{Illustration of the curve $s \mapsto \varphi(\gamma(s),t_0), \; \gamma(0) = x$ for a fixed time $t = t_0$ with vector field \break $s \mapsto \frac{\dif}{\dif s} \varphi(\gamma(s),t) \in T_\xi \Omega_\xi$.}
		\label{yfig2}
	\end{center}
\end{figure}
\begin{multicols}{3}
	\begin{itemize}
		\item $F\colon T_x \Omega_x \to T_\xi \Omega_\xi$
		\item $R\colon T_x \Omega_x \to T_\xi \Omega \xi$
		\item $F^T\colon T_\xi \Omega_\xi \to T_x \Omega_x$
		\item $R^T\colon T_\xi \Omega_\xi \to T_x \Omega_x$
		\item $C = F^T \, F\colon T_x \Omega_x \to T_x \Omega_x$
		\item $B = F \, F^T\colon T_\xi \Omega_\xi \to T_\xi \Omega_\xi$
		\item $\sigma\colon T_\xi \Omega_\xi \to T_\xi \Omega_\xi$
		\item $\tau\colon T_\xi \Omega_\xi \to T_\xi \Omega_\xi$
		\item $S_2 \colon T_x \Omega_x \to T_x \Omega_x$
		\item $S_1 \colon T_x \Omega_x \to T_\xi \Omega_\xi$
		\item $R^T \, \sigma \, R\colon T_x \Omega_x \to T_x \Omega_x$
	\end{itemize}
\end{multicols}

\noindent \textbf{Primary matrix functions} \\
\\
We define primary matrix functions as those functions $\Sigma \colon \Sym^{++}(3) \to \Sym(3)$, such that
\begin{align}
	\Sigma(V) = \Sigma(Q^T \, \text{diag}_V(\lambda_1, \lambda_2, \lambda_3) \, Q) = Q^T \Sigma(\text{diag}_V(\lambda_1, \lambda_2, \lambda_3)) \, Q = Q^T \,
	\begin{pmatrix}
		f(\lambda_1) & 0 & 0 \\
		0 & f(\lambda_2) & 0 \\
		0 & 0 & f(\lambda_3)
	\end{pmatrix} \, Q
\end{align}
with one given real-valued scale-function $f \colon \mathbb{R}_+ \to \mathbb{R}$. Any primary matrix function is an isotropic matrix function but not vice-versa as shows e.g.~$\Sigma(V) = \det V \, \id$. \\
\\
\textbf{List of additional definitions and useful identities}
\begin{itemize}
	\item For two metric spaces $X, Y$ and a linear map $L: X \to Y$ with argument $v \in X$ we write $L.v:=L(v)$. This applies to a second order tensor $A$ and a vector $v$ as $A.v$ as well as a fourth order tensor $\C$ and a second order tensor $H$ as $\C.H$.
	\item We define $J = \det{F}$ and denote by $\Cof X = (\det X)X^{-T}$ the \emph{cofactor} of a matrix in ${\rm GL}^{+}(3)$.
	\item We define $\sym X = \frac12 \, (X + X^T)$ and $\sk X = \frac12 \, (X - X^T)$ as well as $\dev X = X - \frac13 \, \tr(X) \, \id$.
	\item For all vectors $\xi,\eta\in\R^3$ we have the tensor or dyadic product $(\xi\otimes\eta)_{ij}=\xi_i\,\eta_j$.
	\item $S_1=\DD_F \WW(F) = \sigma \, \Cof F$ is the non-symmetric first Piola-Kirchhoff stress tensor.
	\item $S_2=F^{-1}S_1=2\,\DD_C \widetilde{\WW}(C)$ is the symmetric second  Piola-Kirchhoff stress tensor.
	\item $\sigma=\frac{1}{J}\,  S_1\, F^T=\frac{1}{J}\,  F\,S_2\, F^T=\frac{2}{J}\DD_B \widetilde{\WW}(B)\, B=\frac{1}{J}\DD_V \widetilde{\WW}(V)\, V = \frac{1}{J} \, \DD_{\log V} \widehat \WW(\log V)$ is the symmetric Cauchy stress tensor.
	\item $\sigma = \frac{1}{J} \, F\, S_2 \, F^T = \frac{2}{J} \, F \, \DD_C \widetilde{\WW}(C) \, F^T$ is the ''\emph{Doyle-Ericksen formula}'' \cite{doyle1956nonlinear}.
	\item For $\sigma: \Sym(3) \to \Sym(3)$ we denote by $\DD_B \sigma(B)$ with $\sigma(B+H) = \sigma(B) + \DD_B \sigma(B).H + o(H)$ the Fréchet-derivative. For $\sigma: \Sym^+(3) \subset \Sym(3) \to \Sym(3)$ the same applies. Similarly, for $\WW : \R^{3 \times 3} \to \R$ we have $\WW(X + H) = \WW(X) + \langle \DD_X \WW(X), H \rangle + o(H)$.
	\item $\tau = J \, \sigma = 2\, \DD_B \widetilde{\WW}(B)\, B $ is the symmetric Kirchhoff stress tensor.
	\item $\tau = \DD_{\log V} \widehat{\WW}(\log V)$ is the ``\emph{Richter-formula}'' \cite{richter1948,richter1949hauptaufsatze}.
	\item $\sigma_i =\dd\frac{1}{\lambda_1\lambda_2\lambda_3}\dd\lambda_i\frac{\partial g(\lambda_1,\lambda_2,\lambda_3)}{\partial \lambda_i}=\dd\frac{1}{\lambda_j\lambda_k}\dd\frac{\partial g(\lambda_1,\lambda_2,\lambda_3)}{\partial \lambda_i}, \ \ i\neq j\neq k \neq i$ are the principal Cauchy stresses (the eigenvalues of the Cauchy stress tensor $\sigma$), where $g:\mathbb{R}_+^3\to \mathbb{R}$ is the unique function  of the singular values of $U$ (the principal stretches) such that $\WW(F)=\widetilde{\WW}(U)=g(\lambda_1,\lambda_2,\lambda_3)$.
	\item $\sigma_i =\dd\frac{1}{\lambda_1\lambda_2\lambda_3}\frac{\partial \widehat{g}(\log \lambda_1,\log \lambda_2,\log \lambda_3)}{\partial \log \lambda_i}$, where $\widehat{g}:\mathbb{R}^3\to \mathbb{R}$ is the unique function such that \\ \hspace*{0.3cm} $\widehat{g}(\log \lambda_1,\log \lambda_2,\log \lambda_3):=g(\lambda_1,\lambda_2,\lambda_3)$.
	\item $\tau_i =J\, \sigma_i=\dd\lambda_i\frac{\partial g(\lambda_1,\lambda_2,\lambda_3)}{\partial \lambda_i}=\frac{\partial \widehat{g}(\log \lambda_1,\log \lambda_2,\log \lambda_3)}{\partial \log \lambda_i}$ \, 
	\item $T_{\Biot} = \DD_U \widetilde{\WW}(U)$ is the symmetric Biot stress tensor
	\item $\sigma = \frac{1}{\det V} \, V (R \, T_{\Biot}\, R^T) $
	\item $T_{\Biot}^{i} = \frac{\partial g(\lambda_1,\lambda_2,\lambda_3)}{\partial \lambda_i}$ (in case of hyperelasticity)
	\item $\sigma_i = \frac{1}{\lambda_j \, \lambda_k} T_{\Biot}^{i} \; \left(= \; \lambda_i T_{\Biot}^i \; \text{for incompressibility} \right) $
\end{itemize}

\vspace*{2em}
\noindent \textbf{Conventions for fourth-order symmetric operators, minor and major symmetry} \\
\\
For a fourth order linear mapping $\C : \Sym(3) \to \Sym(3)$ we agree on the following convention. \\
\\
We say that $\C$ has \emph{minor symmetry} if
\begin{align}
	\C.S \in \Sym(3) \qquad \forall \, S \in \Sym(3).
\end{align}
This can also be written in index notation as $C_{ijkm} = C_{jikm} = C_{ijmk}$. If we consider a more general fourth order tensor $\C : \R^{3 \times 3} \to \R^{3 \times 3}$ then $\C$ can be transformed having minor symmetry by considering the mapping $X \mapsto \sym(\C. \sym X)$ such that $\C: \R^{3 \times 3} \to \R^{3 \times 3}$ is minor symmetric, if and only if $\C.X = \sym(\C.\sym X)$. \\
\\
We say that $\C$ has \emph{major symmetry} (or is \emph{self-adjoint}, respectively) if
\begin{align}
	\langle \C. S_1, S_2 \rangle = \langle \C. S_2, S_1 \rangle \qquad \forall \, S_1, S_2 \in \Sym(3).
\end{align}
Major symmetry in index notation is understood as $C_{ijkm} = C_{kmij}$. \\
\\
The set of positive-definite, major symmetric fourth order tensors mapping $\R^{3 \times 3} \to \R^{3 \times 3}$ is denoted as $\Sym^{++}_4(9)$, in case of additional minor symmetry, i.e.~mapping $\Sym(3) \to \Sym(3)$ as $\Sym^{++}_4(6)$. By identifying $\Sym(3) \cong \R^6$, we can view $\C$ as a linear mapping in matrix form $\widetilde \C: \R^6 \to \R^6$. \newline If $H \in \Sym(3) \cong \R^6$ has the entries $H_{ij}$, we can write
\begin{align}
	\label{eqvec1}
	h = \text{vec}(H) = (H_{11}, H_{22}, H_{33}, H_{12}, H_{23}, H_{31}) \in \R^6 \qquad \text{so that} \qquad \langle \C.H, H \rangle_{\Sym(3)} = \langle \widetilde \C.h, h \rangle_{\R^6}.
\end{align}
If $\C: \Sym(3) \to \Sym(3)$, we can define $\bfsym \C$ by
\begin{align}
	\langle \C.H, H \rangle_{\Sym(3)} = \langle \widetilde \C.h, h \rangle_{\R^6} = \langle \sym  \, \widetilde \C. h, h \rangle_{\R^6} =: \langle \bfsym \C.H, H \rangle_{\Sym(3)}, \qquad \forall \, H \in \Sym(3).
\end{align}
Major symmetry in these terms can be expressed as $\widetilde \C \in \Sym(6)$. \emph{In this text, however, we omit the tilde-operation and ${\bf sym}$ and write in short $\sym\C\in {\rm Sym}_4(6)$ if no confusion can arise.} In the same manner we speak about $\det \C$ meaning $\det \widetilde \C$. \\
\\
A linear mapping $\C : \R^{3 \times 3} \to \R^{3 \times 3}$ is positive definite if and only if
\begin{align}
	\label{eqposdef1}
	\langle \C.H, H \rangle > 0 \qquad \forall \, H \in \R^{3 \times 3} \qquad \iff \qquad \C \in \Sym^{++}_4(9)
\end{align}
and analogously it is positive semi-definite if and only if
\begin{align}
	\label{eqpossemidef1}
	\langle \C.H, H \rangle \ge 0 \qquad \forall \, H \in \R^{3 \times 3} \qquad \iff \qquad \C \in \Sym^+_4(9).
\end{align}
For $\C: \Sym(3) \to \Sym(3)$, after identifying $\Sym(3) \cong \R^6$, we can reformulate \eqref{eqposdef1} as $\widetilde \C \in \Sym^{++}(6)$ and \eqref{eqpossemidef1} as $\widetilde \C \in \Sym^+(6)$.
\end{appendix}

\end{document}